\documentclass[oneside]{article}
\usepackage[utf8]{inputenc}
\usepackage[english]{babel}  
 \pdfoutput=1 

\usepackage{amsthm}
\newtheorem{theorem}{Theorem}[section]
\newtheorem{lemma}{Lemma}[section]

\newenvironment{remark}[1][Remark]{\begin{trivlist}
		\item[\hskip \labelsep {\bfseries #1}]}{\end{trivlist}}

\theoremstyle{definition}
\newtheorem{definition}{Definition}[section]
\usepackage{amsmath}
\usepackage{amsfonts}
\usepackage{amssymb}
\usepackage{graphicx}
\numberwithin{equation}{section}
\usepackage{lipsum}
\usepackage{mathtools}
\usepackage{scalefnt}
\usepackage{empheq}
\usepackage[a4paper,left=3.5cm,right=3cm,top=3cm,bottom=3cm]{geometry}
\usepackage{makeidx}
\usepackage{bm}

\makeindex

\usepackage{indentfirst}

\begin{document}
	
	\title{Existence and uniqueness for renormalized solutions to a general
		 noncoercive nonlinear parabolic equation}
	
	\author{T.T. Dang\thanks{Department of Mathematics, Leopold Franzens Universit\"at Innsbruck, Austria, e-mail: Thi-Tam.Dang@uibk.ac.at}, G. Orlandi\thanks{Dipartimento di Informatica, Universit\`a di Verona, Italy, e-mail: giandomenico.orlandi@univr.it} \\}
	
	\maketitle
	
\begin{abstract}
This paper introduces the concept of renormalized solution for a general class of noncoercive nonlinear parabolic problems, including both singularities and unbounded lower order terms. We prove existence and uniqueness of renormalized solutions for such class of problems.  

\end{abstract}

\section{Introduction}
Let $\Omega$ be an open bounded subset of $\mathbb{R}^{N}$ with $N\geq 2$, for any $T>0$, and denote $Q := \Omega \times (0,T)$ and $\Gamma_{T} := \partial \Omega \times (0,T)$. This paper is concerned with nonlinear parabolic problems of the form
\begin{align}\label{1.1}
\begin{cases}
	&	u_{t} - \mathrm{div}A(x, t, u,\nabla u)  + H(x,t,u,\nabla u) + G(x,t,u) = - \mathrm{div}(F(x,t)) \quad \text{in} \ \  Q, \\
	&	u = 0  \quad \text{on} \quad \Gamma_{T},\\ 
	& u(\cdot, 0) = u_0 \quad \text{in } \Omega,
\end{cases}	
\end{align}
where the initial datum $u_0 \in L^2(\Omega)$, the vector field $F(x,t) \in L^{p^{\prime}}(0,T;W^{-1, p^{\prime}}(\Omega))$ for any $2^*:=\frac{2N}{N+2} <p <N$ and the operator $A= A(x,t, z,\xi): Q \times \mathbb{R} \times\mathbb{R}^{N} \rightarrow \mathbb{R}^{N}$ is a Carath\'eodory function satisfying, for a.e. $(x,t)\in Q$, $z \in \mathbb{R}$ and for any $\xi, \xi^{*}  \in \mathbb{R}^{N}$, the following growth conditions:
\begin{align}
	&	\left\langle A(x,t, z, \xi) \xi  \right\rangle \geq \alpha \vert \xi \vert^{p}- \Big( b(x,t)\vert z \vert \Big)^{p}- K(x,t), \quad \alpha >0, \label{1.2}\\
	& \vert A(x,t,z, \xi) \vert \le \beta | \xi |^{p-1} + \Big(b(x,t)\vert z \vert\Big)^{p-1}+ a_0(x,t), \quad \beta >0, \label{1.3}\\
	& 	\left\langle A(x,t,z,\xi)- A(x,t,z,\xi^{*}), \xi - \xi^{*} \right\rangle  \ge 0, \qquad \xi \ne \xi^{*}, \label{1.4}
\end{align}
where the term $b(x,t) \in L^{\infty}(0,T; L^{N,\infty}(\Omega))$, $K(x,t) \in L^1(
Q)$ and $a_0 \in L^{p^{\prime}} (Q)$. Here $L^{N,q}(\Omega)$, for $1\le q\le \infty$, denotes the Lorentz space (see Sect. \ref{sect:lorentz} for definitions).

Moreover, for a.e. $(x,t)\in Q$ and for every $z \in \mathbb{R}$ and $\xi \in \mathbb{R}^{N}$, the term $H= H(x,t,z, \xi): Q \times \mathbb{R} \times \mathbb{R}^{N} \rightarrow \mathbb{R} $ is a Carath\'eodory function satisfying 
\begin{align}\label{1.5}
	\vert H(x,t,z,\xi) \vert \le c_0(x,t) \vert \xi \vert^{\gamma} + c_1(x,t), \quad 0 \le \gamma <p-1,
\end{align}
where $	c_0(x,t) \in L^{\infty} (0,T;L^{N,1}(\Omega))$ and $c_1(x,t) \in L^1(Q)$.

Finally, the term $G= G(x,t,z): Q \times \mathbb{R} \rightarrow \mathbb{R}$ is a Carath\'eodory function that satisfies the following conditions:
\begin{align}
	& \vert G(x,t,z) \vert \le d_0(x,t) \vert z \vert^{\lambda} + d_1(x,t), \quad 	0 \le \lambda < \frac{N(p-1)}{p-1}, \label{1.6}\\	
	&G(x,t,z)z \ge 0, \label{1.7}
\end{align}
 where $d_0 \in L^{\infty} (0,T, L^{r^{\prime},1}(\Omega))$ with $r =  \frac{N(p-1)}{\lambda (n-p)}$ and $d_1(x,t) \in L^1(Q)$.
 

In case $p=2$, $H=G=0$, and classical growth conditions on the operator $A$, problem \eqref{1.1} corresponds to a homogeneous Fokker-Planck equation (see e.g. \cite{A. Porretta}), that characterizes the evolution of certain Brownian motion processes.

In related studies by Farroni and Greco $\cite{Farroni-Greco}$ and Cardaliaguet {\it et al. }$\cite{Cardaliaguet}$, problem \eqref{1.1} is considered in the case $H=G=0$ and also assuming respectively $b(t,x)=0$ or $b(x,t) \in L^{\infty}(0,T; L^1(\Omega))$ for the growth conditions \eqref{1.2},\eqref{1.3} on the operator $A$. In their works, existence of solutions is proved first for a sequence of approximating  problems, then by showing that limits of approximate solutions converge to solutions of the original problem. 

In this paper, we consider a more general class of nonlinear evolution problems, including nonzero terms $H$ and $G$ in \eqref{1.1}, that play the role of respectively a convection and a reaction term, and a more general unbounded term $b(t,x)\in L^{\infty}(0,T; L^{N,\infty}(\Omega))$ in \eqref{1.2},\eqref{1.3}. This generalization presents a twofold complexity. First, the possible singular behavior of the operator $A$, arising from the specific assumptions made on the function $b(x,t)$. This singular behavior poses significant obstacles in the formulation and analysis of the problem and requires special techniques to deal with its effects. Second, the possibly unbounded convective coefficient $c_0(x,t) \in L^{\infty} (0,T; L^{N,1}(\Omega))$ in \eqref{1.5} adds further complexity to the problem.

Moreover, no coerciveness assumption is imposed on the operator $A$, thus it is not feasible to obtain a distributional solution to the problem \eqref{1.1}. To overcome this challenge, we employ the framework of renormalized solutions within the context of this study. The concept of renormalized solutions was initially introduced by Lions and Di Perna (see $\cite{Lions}$) in their investigation of the Boltzmann equation. Subsequently, it was adapted for nonlinear elliptic problems (see $\cite{Betta}$). Simultaneously, the notion of entropy solution was independently developed by Benilan (see, $\cite{Benilan}$) for the study of elliptic problems. These theoretical frameworks provide alternative approaches to address the nonexistence of distributional solutions, thus enabling a comprehensive analysis of the problem at hand.

Regarding the framework of renormalized solutions in the parabolic case, existence and uniqueness results have been established under specific conditions. In particular, when the operator $A(x,t,z,\xi)$ does not depend on the variable $z$ and the functions $H,G$ are identically zero, existence and uniqueness of renormalized solutions have been proven in $\cite{Blanchard, Murat, Boccardo}$). An existence result for a similar problem involving a nonlinear term with a growth condition proportional to $\vert Du \vert^{p}$ was proved by Porretta $\cite{Porettea A}$. In the recent work $\cite{Dang}$ it is shown existence and uniqueness of renormalized solutions for a noncoercive parabolic problem of convection-diffusion type. More precisely, the study  in \cite{Dang} focuses on the scenario where the term $A(x,t,z,\xi)$ has no dependence on the variable $z$, the convective term involves an unbounded coefficient in the Lorentz class, and the functions $H,G$ are identically zero. 
 
This paper presents further advancements compared to prior research endeavors, offering novel contributions in two distinct key aspects.  Firstly, we delve into the investigation of the growth behavior of the operator $A(x,t,z,\xi)$ with respect to the variable $z$ through the unbounded function $b(t,x)$, which introduces additional complexities and challenges in the analytical investigation. This constitutes a novel aspect within the context of our study, setting it apart from previous works. Secondly, we address the presence of a possibly unbounded term in the growth condition \eqref{1.5} for the function $H$. Our main results, rigorously stated in Theorem \ref{Th: 3.1} in Section 3 and Theorem \ref{Th: 4.1} in Section 4, establish respectively existence and uniqueness of renormalized solutions to problem \eqref{1.1} imposing further assumptions on the  unbounded coefficient $c_0(t,x)$, namely condition \eqref{3.7}. This enables us to overcome the challenges associated with the unboundedness of the coefficient, and provides a suitable framework for analyzing the problem under consideration. The main difficulty encountered in our proof lies in establishing an a priori estimate for the quantity $\| \vert \nabla u \vert^{p-1} \|_{L^{N^{\prime},\infty}(\Omega)}$.

The structure of the paper is outlined as follows. Section~2 provides an overview of Lorentz spaces, which form an essential component of the mathematical framework employed in this study. In Section~3 we focus on establishing the existence of renormalized solutions to problem \eqref{1.1}. The proof of our main existence result, Theorem \ref{Th: 3.1}, is accomplished through the application of a series of technical lemmas. Section~4 investigates the uniqueness of the renormalized solution, stated in Theorem \ref{Th: 4.1}. A rigorous proof of uniqueness is presented as an immediate consequence of a suitable comparison lemma (Lemma \ref{Le: 4.1}). In the Appendix we gather the technical details to complete the proof of Lemma \ref{Le: 4.1}.

\section{Lorentz spaces}\label{sect:lorentz}
In this section, we recall the definition and some properties of Lorentz spaces, which are used extensively in our work. For a detailed background of Lorentz spaces, we refer the reader to $\cite{Lorentz}$, where in-depth discussions and references can be found.

We denote the distribution function of $f$ by $\omega_{f}(h) = \mathrm{meas} \left\lbrace  x \in \Omega: \vert f(x) \vert >h  \right\rbrace$. The Lorentz space $L^{p,q}(\Omega)$, for $1\le p,q<+\infty,$ consists of all measurable functions $f$ defined on $\Omega$ such that
\begin{align}
	\Vert f \Vert_{p,q} = \left(p \int_0^{\infty} [\omega_{f}(h)]^{\frac{p}{q}} h^{q-1} dh \right)^{\frac{1}{q}} \, < +\infty\, . 
\end{align}
The Lorentz space $L^{p,q}(\Omega)$, endowed with the norm $\| \cdot \|_{p,q}$ is a Banach space. In case $p=q$, $L^{p,q}(\Omega)$ coincides with the Lebesgue space $L^{p}(\Omega)$. For $q = \infty$, the class $L^{p,\infty}(\Omega)$ consists of all measurable functions $f$ defined on $\Omega$ for which
\begin{align}
	\Vert f \Vert _{p, \infty}^{p} = \sup_{h>0} h^{p} \omega_{f}(h) < + \infty
\end{align}
and it coincides with the Marcinkiewicz class weak-$L^{p}(\Omega)$.

 For every $1 \le q<p<r \le +\infty$, the following continuous inclusions hold true:
\begin{align}\label{2.3}
	L^{r}(\Omega) \subset L^{p,q}(\Omega) \subset L^{p,r}(\Omega) \subset L^{p,\infty}(\Omega) \subset L^{q}(\Omega).
\end{align}
Recall the generalized H\"older inequality valid for Lorentz spaces: for $1<p_1,p_2<+\infty$, $1\leq q_1,q_2\le +\infty$,  $f \in L^{p_1, q_1}(\Omega)$ and $g \in L^{p_2, q_2}(\Omega)$, 
 we have
\begin{align}
	\| f g\|_{p,q} \le \| f \|_{p_1, q_1} \| g \|_{p_2, q_2},
\end{align}
where $\frac{1}{p} = \frac{1}{p_1}+\frac{1}{p_2}$ and $\frac{1}{q} = \frac{1}{q_1}+\frac{1}{q_2}$.
Denote by
\begin{align}
	\mathrm{dist}_{L^{p, \infty}(\Omega)} (f, L^{\infty}(\Omega)) = \inf_{g \in L^{\infty}(\Omega)} \| f-g \|_{p, \infty}
\end{align}
the distance of a function $f \in L^{p,\infty}$ to $L^{\infty}(\Omega)$. We have the  characterization
\begin{align}
 \mathrm{dist}_{L^{p, \infty}(\Omega)} (f, L^{\infty}(\Omega)) =\lim_{m \rightarrow +\infty} \| R_{m}f\|_{p, \infty},
\end{align}
where, for $f\in L^{p, \infty}(\Omega)$, we define the operator $R_m: L^{p, \infty}(\Omega)\to L^{p, \infty}(\Omega)$ as
\begin{equation}\label{eq: remainder}
R_m(f)=f- T_{m}f\, ,
\end{equation}
with $T_{m}$ the truncation operator at level $m>0$, i.e. $T_m(f)(x)=\sup \{ -m, \inf\{ f(x), m\} \}$.

The generalization of Sobolev embedding theorem to Lorentz spaces reads as follows: 
	\begin{theorem}$\cite{Farroni,Lorentz,O'Neil}$\label{Th: 2.1}
	Assume that $1<p<N, 1\le q \le p$. Then, any function $g\in W_0^{1,1}(\Omega)$ satisfying $\vert\nabla g \vert \in L^{p,q}(\Omega)$ belongs to $L^{p^{\star},q}(\Omega)$ where $p^{*}= \frac{Np}{N-p}$ and
	\begin{align}
		\Vert g \Vert_{p^{*},q} \le S_{N,p} \Vert \nabla g\Vert_{p,q}
	\end{align}
	where $S_{N,p}= \omega_{N}^{-1/N} \frac{p}{N-p}$ and $\omega_{N}$ stands for the measure of the unit ball in $\mathbb{R}^{N}$.
\end{theorem}
	
\section{Existence of renormalized solutions}	
In this section we give the definition of renormalized solution for problem \eqref{1.1}. Moreover, we will prove existence of a renormalized solution assuming suitable bounds on the lower order terms of the nonlinear parabolic equation in \eqref{1.1}.
\begin{definition}(Renormalized solution)\label{Def: 3.1}
	A real-valued function $u$ defined on $Q$ is a renormalized solution of problem \eqref{1.1} if for any $k,n>0$,
	\begin{align}
		&u\in C^0([0,T],L^2(\Omega)) \cap L^{p}([0,T]; W_0^{1,p}(\Omega)),\\
		&T_{k}(u)\in L^{p}(0,T; W_0^{1,p}(\Omega)), 
	\end{align}	
	\begin{align}
		T_{k+n}(u)- T_{k}(u) \rightarrow 0 \ \ \text{strongly in} \ \ L^{p}(0,T; W_0^{1,p}(\Omega))\quad\text{as } k\to+\infty,
	\end{align}
	and for any $S\in C^{\infty}(\mathbb{R})$ such that $S^{\prime} \in C_0^{\infty}(\mathbb{R})$, we have
	\begin{equation}\label{3.4}
		\begin{aligned}
		&	(S(u))_{t}  - \mathrm{div} \Big(  A(x,t,u, \nabla u)\cdot S^{\prime}(u) \Big)
			 + S^{\prime \prime}(u) A(x,t,u,\nabla u) \cdot \nabla u  	\\
		&\qquad	+H(x,t.u, \nabla u)S^{\prime}(u) + G(x,t, u) S^{\prime}(u) = \mathrm{div} \Big( F  S^{\prime}(u) \Big) - S^{\prime \prime}(u) F \cdot \nabla u,\\
		\end{aligned}
	\end{equation}
  and
	\begin{align} \label{3.5}
		S(u(t=0))=S(u_0).
	\end{align}
\end{definition}

 Let us further define, for any $ \theta \in [0,1]$,  
\begin{align}\label{3.6}
	\vartheta_{\theta} (x,t)= \theta c_0(x,t) + (1- \theta) b(x,t). \ 
\end{align}
	

Our main existence result is stated in the following theorem.
\begin{theorem}[Existence]\label{Th: 3.1}
	Assume assumptions \eqref{1.2}-\eqref{1.7} are fulfilled. Furthermore, assume that, for some $m>0$,
\begin{equation}\label{3.7}
	\begin{aligned}
		\| 	R_{m}(	\vartheta_{\theta}) \|_{L^\infty (0,T; L^{N,\infty}(\Omega))} \le \frac{ 1 }{S_{N,p}}\left(\frac{\alpha}{2p} \right)^{1/p} , 
	\end{aligned}
\end{equation}
where $\vartheta_\theta$ is defined in \eqref{3.6} and $R_m$ is defined in \eqref{eq: remainder}. Then  problem \eqref{1.1} has at least one renormalized solution according to Definition \ref{Def: 3.1}.
\end{theorem}
\begin{remark}\label{rem: 3.1} Condition \eqref{3.7} implies in particular, by \eqref{2.3}, 
\begin{equation}\label{3.8}
	\begin{aligned}
		\| R_{m}(	\vartheta_{\theta}) \|_{L^\infty (0,T; L^{N,q}(\Omega))} \le \frac{ 1 }{S_{N,p}}\left(\frac{\alpha}{2p} \right)^{1/p}  \text{for any } 1\le q \le +\infty.
	\end{aligned}
\end{equation}
\end{remark}

\subsection{The approximating problems}
The strategy for constructing a renormalized solution follows the ideas in $\cite{Blanchard, Murat, Porettea A}$. We consider the following approximating problems for \eqref{1.1} :
\begin{align} \label{3.8}
	\begin{cases}
		(u_{\epsilon})_{t}- \mathrm{div} A_{\epsilon}(x, t,u_{\epsilon}, \nabla u_{\epsilon}) + H_{\epsilon}(x,t,u_{\epsilon}, \nabla u_{\epsilon})+ G_{\epsilon}(x,t, u_{\epsilon}) = - \mathrm{div} F_{\epsilon} \quad   \ \text{in}  \quad Q, \\
		u_{\epsilon} = 0, \quad \text{on} \ \ \Gamma_{T},\\
		u_{\epsilon}(\cdot, 0)= u_{0 \epsilon} \quad \text{in} \quad \Omega,
	\end{cases}
\end{align}
where
\begin{equation}
	\begin{aligned}
	   & A_{\epsilon}(x,t,z,\xi) = A(x, t, T_{1/\epsilon}(z), \xi),\\
		&F_{\epsilon}= T_{1/ \epsilon}(F),\\
		& u_{0 \epsilon} = T_{1/ \epsilon}(u_0),\\
		& H_{\epsilon} (x,t,z, \xi) = T_{1/ \epsilon} H(x,t,z ,\xi),\\
		& G_{\epsilon}(x,t, z)  = T_{1/ \epsilon} G(x,t,z).
	\end{aligned}
\end{equation}
The operator $A_{\epsilon}(x,t,u_{\epsilon}, \nabla u_{\epsilon})$ satisfies assumptions \eqref{1.2}-\eqref{1.4}, 
that we express here in term of $\vartheta_\theta$:
\begin{align}
	& \left\langle A_{\epsilon}(x,t,z, \xi), \xi\right\rangle  \ge \alpha \vert \xi \vert^{p} - \Big( \vartheta_0(x,t) \vert z \vert \Big)^{p}- K(x,t), \label{3.10}\\
 	& \vert A_{\epsilon}(x,t,z, \xi) \le \beta | \xi |^{p-1} + \left(  \vartheta_0(x,t) \vert z \vert  \right) ^{p-1}+ a_0(x,t),  \label{3.11}\\
	& 	\left\langle A_{\epsilon}(x,t,z,\xi)- A_{\epsilon}(x,t,z,\xi^{*}), \xi - \xi^{*} \right\rangle  \ge 0, \qquad \xi \ne \xi^*,\label{3.12}
\end{align}
where $ \vartheta_0(x,t) \in L^{\infty}(0,T; L^{N, \infty}(\Omega))$ is given by \eqref{3.6}.

By assumption \eqref{1.5} on $H$ the approximating term $H_{\epsilon}$ satisfies the bounds
\begin{align}
	&	\vert H_{\epsilon}(x,t,z, \xi) \vert \le \vert H(x,t,z,\xi)\vert \le \vartheta_1(x,t) \vert \xi \vert^{\gamma}+c_1(x,t), \label{3.13}\\
& \vert  H_{\epsilon}(x,t,z, \xi) \vert \le \frac{1}{\epsilon},	\label{3.14}
\end{align}
with $  \vartheta_1(x,t) \in L^{\infty}(0,T; L^{N,1}(\Omega) ) $.

Furthermore, assumptions \eqref{1.6}-\eqref{1.7} on $G$ yield the following conditions: 
\begin{align}
	& G_{\epsilon}(x,t,z)z \ge 0, \label{3.15}\\
	& \vert G_{\epsilon}(x,t,z) \vert \le \vert G(x,t,z) \vert \le d_1(x,t) \vert z \vert^{\lambda} +d_2(x,t) \label{3.16}\\
	& \vert G_{\epsilon}(x,t,z) \vert \le \frac{1}{\epsilon} \label{3.17}.
\end{align}
Since $F_{\epsilon}\in L^{p^{\prime}} (0,T; W^{-1, p^{\prime} }(\Omega) )$ and $u_{0\epsilon}\in L^2(\Omega)$, by \eqref{3.14}, \eqref{3.15} and \eqref{3.17} the approximating problem \eqref{3.8} admits a solution $u_{\epsilon}\in L^{p}(0,T; W_0^{1,p}(\Omega)) \cap C^0(0,T,L^2(\Omega))$, see Theorem~2.1 in \cite{Porettea A} for a detailed proof.

\subsection{A priori estimates}

To prove Theorem \ref{Th: 3.1}, we first give an estimate on $ R_{k}(u_{\epsilon})$, where $k$ is a positive number depending on $\epsilon$ and the data (see Lemma \ref{Le: 3.3}).
Taking advantage of Lemma \ref{Le: 3.3} we are able to provide bounds on the truncation $T_{k}(u_{\epsilon})$ according to Lemma \ref{Le: 3.2} below, which extends the results of $\cite{Betta}$ to the context of noncoercive nonlinear parabolic problems. The subsequent lemmas of this section provide compactness arguments to show the existence of a convergent subsequence of approximate solutions, and preparatory results to show that limits of approximate solutions are indeed renormalized solutions to \eqref{1.1}.

\begin{lemma}\label{Le: 3.3}
	Let the assumptions in Theorem \ref{Th: 3.1} be satisfied. Then every solution $u_{\epsilon}$ of \eqref{3.8} satisfies
	\begin{align}
		\| \vert \nabla R_{k}(u_{\epsilon})\vert^{p-1} \|_{L^1(0,T; L^{N^{\prime}, \infty}(\Omega) )} \le C^{\prime}(N,p) \left[ \tilde{M}+ \vert \Omega \vert^{\frac{1}{N^{\prime}}- \frac{1}{p^{\prime}} } \tilde{L}^{\frac{1}{p^{\prime}}} \right] ,
	\end{align}
where $\tilde{M}, \tilde{L}$ are given constants. Further, for any $n \ge 1$, the following statement holds true:
	\begin{align}\label{3.25}
		\lim_{n\rightarrow \infty} \int_0^{T}\int_{\left\lbrace n\le \vert u_{\epsilon}\vert \le n+1 \right\rbrace } \vert \nabla u_{\epsilon}\vert^{p} dxdt=0.
	\end{align}
\end{lemma}

\begin{proof}
	For any $n \ge 1$, we define the Lipschitz continuous function $\phi_{n}$ as follows:
	\begin{align}\label{phi}
		\phi_{n}(r) = T_{n+1}(r) - T_{n}(r).
	\end{align}	
	Using the test function $ \phi_{n}(u_{\epsilon})$ in the approximating problem \eqref{3.8} and setting 
	\begin{align}\label{tildephi}
		\tilde{\phi}_{n}(s) = \int_0^{s} \phi_{n}(\xi)d\xi ,
	\end{align} 
	we have
	\begin{equation}\label{3.26}
		\begin{aligned}
			&	\int_{\Omega} \tilde{\phi}_{n}(u_{\epsilon}(x,T))dx + \int_{Q} A_{\epsilon}(x,t,u_{\epsilon},\nabla u_{\epsilon}) \cdot \nabla \phi_{n}(u_{\epsilon})dx dt\\
			& \quad+  \int_{Q} H_{\epsilon} (x,t,u_{\epsilon}, \nabla u_{\epsilon}) \phi_{n}(u_{\epsilon}) dxdt +  \int_{Q} G_{\epsilon}(x,t,u_{\epsilon}) \phi_{n}(u_{\epsilon}) dxdt\\
			& = \int_{\Omega} \tilde{\phi}_{n}(u_{0 \epsilon}) dx+    \int_{Q} F_{\epsilon} \cdot \nabla \phi_{n}(u_{\epsilon}) dxdt .
		\end{aligned}
	\end{equation}


	Observe that $\nabla \phi_{n}(u_{\epsilon}) = \chi_{\left\lbrace n \le \vert u_{\epsilon} \vert \le n+1 \right\rbrace } \nabla u_{\epsilon}$. By the coercivity assumption \eqref{3.10} on $A_{\epsilon}$ we obtain
	\begin{equation}\label{3.27}
		\begin{aligned}
			&\int_{Q} A_{\epsilon}(x,t,u_{\epsilon},\nabla u_{\epsilon}) \cdot \nabla \phi_{n}(u_{\epsilon})dxdt\\
			& \quad = \int_0^{T}\int_{ \left\lbrace n \le \vert u_{\epsilon} \vert \le n+1 \right\rbrace }  A_{\epsilon}(x,t,u_{\epsilon}, \nabla u_{\epsilon}) \cdot \nabla u_{\epsilon} dx dt\\
			& \quad \ge \alpha \int_0^{T} \int_{ \left\lbrace n \le \vert u_{\epsilon} \vert \le n+1 \right\rbrace } \vert \nabla u_{\epsilon} \vert^{p}dxdt\\
			& \qquad -  \int_0^{T}\int_{ \left\lbrace n\le \vert u_{\epsilon} \vert \le n+1 \right\rbrace } \Big( \vartheta_0(x,t) \vert u_{\epsilon} \vert \Big)^{p}dx dt -  \int_0^{T} \int_{ \left\lbrace  \vert u_{\epsilon} \vert >n\right\rbrace } \vert K(x,t) \vert dxdt.
		\end{aligned}
	\end{equation}
	Thanks to H\"older inequality and  Sobolev embedding, we further obtain
	\begin{equation}\label{3.28}
		\begin{aligned}
			& \int_0^{T} \int_{ \left\lbrace n \le \vert u_{\epsilon} \vert \le n+1 \right\rbrace } \Big( \vartheta_0(x,t) \vert u_{\epsilon} \vert \Big)^{p}dxdt \\
			& \quad \le \int_0^{T} \int_{ \left\lbrace n< \vert u_{\epsilon} \vert <n+1 \right\rbrace } \Big( (R_{m} \vartheta_0(x,t) ) \vert u_{\epsilon}\vert \Big)^{p} dxdt \\
			&  \qquad+  \int_0^{T} \int_{ \left\lbrace n< \vert u_{\epsilon} \vert <n+1 \right\rbrace }  \Big( T_{m} \vartheta_0(x,t) \vert u_{\epsilon} \vert \Big)^{p} dx dt\\
			& \quad \le \| R_{m} \vartheta_0 \|^{p}_{L^{\infty}(0,T, L^{N,\infty}(\Omega))} \int_0^{T}  \int_{ \left\lbrace n< \vert u_{\epsilon} \vert < n+1 \right\rbrace } \vert u_{\epsilon} \vert^{p}dxdt\\
			&  \qquad+ m^{p} \int_0^{T}  \int_{\left\lbrace  n< \vert u_{\epsilon} \vert < n+1\right\rbrace}	\vert u_{\epsilon} \vert^{p} dx dt\\
			& \quad \le  S_{N,p}^{p} \| R_{m} \vartheta_0 \|^{p}_{L^{\infty}(0,T;L^{N,\infty}(\Omega))} \int_0^{T} \int_{ \left\lbrace n< \vert u_{\epsilon} \vert <n +1 \right\rbrace } \vert \nabla u_{\epsilon} \vert^{p}dxdt\\
			&  \qquad + m^{p}  \| 1 \|^{p}_{L^{N,\infty}(\Omega)} S_{N,p}^{p}  \int_0^{T} \int_{ \left\lbrace n< \vert u_{\epsilon} \vert <n +1 \right\rbrace } \vert \nabla u_{\epsilon} \vert^{p}dxdt.
		\end{aligned}
	\end{equation}
	By assumption \eqref{3.7} we have
	\begin{align*}
		S_{N,p}^{p} \| R_{m} \vartheta_0 \|^{p}_{L^{\infty}(0,T; L^{N,\infty}(\Omega))} \le \frac{\alpha}{2p}.
	\end{align*}
	Therefore
	\begin{equation}\label{3.29}
		\begin{aligned}
			& \int_0^{T} \int_{ \left\lbrace n \le \vert u_{\epsilon} \vert \le n+1 \right\rbrace } \Big( \vartheta_0(x,t) \vert u_{\epsilon} \vert \Big)^{p}dxdt \\
			& \le \Big( \frac{\alpha}{2p}+ m^{p}\vert \Omega \vert^{\frac{p}{N}}S_{N,p}^{p}  \Big) \int_0^{T} \int_{ \left\lbrace n \le \vert u_{\epsilon} \vert \le n+1 \right\rbrace }  \vert \nabla u_{\epsilon} \vert^{p} dxdt.
		\end{aligned}
	\end{equation}
	
	Using $ \phi_{n}(u_{\epsilon})= 0$ for $\vert u_{\epsilon} \vert \le n$, assumption \eqref{3.7}, the growth assumption \eqref{3.13} on $H_{\epsilon}$, and the generalized H\"older inequality, we get
	
	\begin{equation*}\label{3.30}
		\begin{aligned}
			&	\left|  \int_{Q} H_{\epsilon}(x,t,u_{\epsilon},\nabla u_{\epsilon}) \phi_{n}(u_{\epsilon})dxdt \right|\\
			& \quad \le  \int_0^{T}  \int_{ \left\lbrace \vert u_{\epsilon} \vert >n \right\rbrace } \vert H_{\epsilon}(x,t,u_{\epsilon},\nabla u_{\epsilon}) \phi_{n}(u_{\epsilon}) \vert dxdt\\
			& \quad \le  \int_0^{T}  \int_{ \left\lbrace \vert u_{\epsilon} \vert >n \right\rbrace } \Big( \vartheta_1 (x,t)\vert \nabla u_{\epsilon}\vert^{p-1}  +  c_1(x, t) \Big) dx dt \\
			& \quad \le \int_0^{T} \int_{ \left\lbrace \vert u_{\epsilon} \vert >n \right\rbrace } (R_{m} \vartheta_1(x,t)) \cdot \vert \nabla u_{\epsilon} \vert^{p-1} dxdt \\
			&  \qquad+  \int_0^{T}  \int_{ \left\lbrace \vert u_{\epsilon} \vert >n \right\rbrace } T_{m} \vartheta_1(x,t) \cdot \vert \nabla u_{\epsilon}\vert^{p-1}dxdt   + \int_0^{T}  \int_{ \left\lbrace \vert u_{\epsilon} \vert >n \right\rbrace } c_1(x,t) dxdt\\
			& \quad \le  \| R_{m} \vartheta_1\|_{L^{\infty}(0,T; L^{N,1}(\Omega))} \int_0^{T}  \| \vert \nabla R_{n}(u_{\epsilon})\vert^{p-1}\|_{L^{N^{\prime}, \infty}(\Omega)} dt \\
			& \qquad + m \int_0^{T}  \| \vert \nabla R_{n}(u_{\epsilon})\vert^{p-1}\|_{L^{N^{\prime}, \infty}(\Omega)} dt + \int_0^{T} \int_{ \left\lbrace \vert u_{\epsilon} \vert >n \right\rbrace } c_1(x,t)dx dt\\
			& \quad \le \eta \int_0^{T}  \| \vert \nabla R_{n}(u_{\epsilon})\vert^{p-1}\|_{L^{N^{\prime}, \infty}(\Omega)} dt +  \int_0^{T} \int_{ \left\lbrace \vert u_{\epsilon} \vert >n \right\rbrace } c_1(x,t)dx dt,
		\end{aligned}
	\end{equation*}
where $\eta$ is given by
\begin{align}\label{3.23a}
	\eta =  \frac{1}{S_{N,p}}\left(  \frac{\alpha}{2p}\right)^{\frac{1}{p}}+ m.
\end{align}

	The sign condition \eqref{3.15} on $G_{\epsilon}(x,t,u_{\epsilon})$ leads to
	\begin{equation}\label{3.31}
		\begin{aligned}
			\int_{Q} G_{\epsilon} (x,t,u_{\epsilon},\nabla u_{\epsilon}) \phi_{n}(u_{\epsilon})dxdt \ge 0.
		\end{aligned}
	\end{equation}
	
	By virtue of Young's inequality, we have
	\begin{equation}\label{3.32}
		\begin{aligned}
			&	 \int_{Q} F_{\epsilon} \cdot \nabla \phi_{n}(u_{\epsilon})dx dt\\
			& \quad	\le \frac{\alpha}{2p} \int_0^{T} \int_{ \left\lbrace n \le \vert u_{\epsilon} \vert \le n+1 \right\rbrace }  \vert \nabla u_{\epsilon} \vert^{p} dxdt + \frac{2^{p^{\prime}/p}}{p^{\prime}\alpha^{p^{\prime}/p}} \int_0^{T}  \int_{ \left\lbrace \vert u_{\epsilon} \vert > n\right\rbrace }
			\vert F_\epsilon\vert^{p^{\prime}} dxdt.
		\end{aligned}
	\end{equation}	
Denote $C_1 = \frac{\alpha}{p^{\prime}} - m^{p} \vert \Omega \vert^{\frac{p}{N}} S_{N,p}^{p}.$ Combining \eqref{3.27}-\eqref{3.32}, we obtain
	\begin{equation}\label{3.33}
		\begin{aligned}
			&	\int_{\Omega} \tilde{\phi}_{n}(u_{\epsilon}(x,T))dx + C_1 \int_0^{T} \int_{ \left\lbrace n \le \vert u_{\epsilon} \vert \le n+1 \right\rbrace } \vert \nabla u_{\epsilon} \vert^{p} dxdt  \\
			& \le  \int_{\Omega}  \tilde{\phi}_{n}(u_{0 \epsilon})dx +  M + L,
		\end{aligned}
	\end{equation}	
	where
	\begin{equation*}
		\begin{aligned}
			&M = \eta \int_0^{T}   \| \vert \nabla R_{n}(u_{\epsilon})\vert^{p-1}\|_{L^{N^{\prime}, \infty}(\Omega)} dt   +  \int_0^{T} \int_{ \left\lbrace \vert u_{\epsilon} \vert >n \right\rbrace } c_1(x,t)dx dt ,\\
			& L =  \frac{2^{p^{\prime}/p}}{p^{\prime}\alpha^{p^{\prime}/p}} \int_0^{T} \int_{ \left\lbrace \vert u_{\epsilon} \vert >n \right\rbrace }	\vert F_\epsilon\vert^{p^{\prime}} dx dt+ \int_0^{T}\int_{ \left\lbrace \vert u_{\epsilon} \vert >n \right\rbrace } K(x,t) dxdt.
		\end{aligned}
	\end{equation*}
	
	Since $\phi_{n}(u_{\epsilon})$ converges to $0$ pointwise as $n \rightarrow + \infty$,  the sequence $\phi_{n}(u_{\epsilon})$ converges to $0$ weakly in $L^{p}(0,T, W_0^{1,p}(\Omega))$. Moreover $\tilde{\phi}_{n}(u_{0\epsilon}) \rightarrow 0$ a.e. in $\Omega$ as $n \rightarrow + \infty$ and $\vert \tilde{\phi}_{n}(u_{0\epsilon}) \vert \le u_{0 \epsilon}$ a.e. in  $\Omega$. Since $u_{0\epsilon} \in L^2(\Omega)$, Lebesgue's convergence theorem implies that
	\begin{align}
		\int_{\Omega} \tilde{\phi}_{n}(u_{0 \epsilon}) dx \rightarrow 0 
	\end{align}
	as $n \rightarrow + \infty$.
	
	Therefore, we have
	\begin{equation}\label{eq: 3.43}
		\begin{aligned}
			\int_0^{T}\int_{ \left\lbrace n \le \vert u_{\epsilon} \vert \le n+1 \right\rbrace } \vert \nabla u_{\epsilon} \vert^{p} dxdt  \le \tilde{M} + \tilde{L},
		\end{aligned}
	\end{equation}
	where
	\begin{equation}
		\tilde{M} = \frac{1}{C_1} M,\qquad	\tilde{L} = \frac{1}{C_1} L.
	\end{equation}
	From Lemma \ref{Le: 3.2} below it follows that
	\begin{equation}\label{eq: 3.45}
		\begin{aligned}
			& \| \vert \nabla R_{n}(u_{\epsilon})\vert^{p-1}\|_{L^1(0,T;L^{N^{\prime},\infty }(\Omega)) }  \\
			&  \quad \le C^{\prime}(N,p)\left[ \tilde{M} + \vert\Omega \vert^{\frac{1}{N^{\prime}}- \frac{1}{p^{\prime}} } \tilde{L}^{\frac{1}{p^{\prime}}} \right]\\
			& \quad \le C^{\prime}(N,p) \left[  \frac{\eta}{C_1} \| \vert \nabla R_{n}(u_{\epsilon})\vert^{p-1}\|_{L^1(0,T; L^{N^{\prime},\infty }(\Omega)) }   \right] \\
			&   \qquad+ C^{\prime}(N,p) \left[  \frac{1}{C_1} \int_0^{T} \int_{ \left\lbrace \vert u_{\epsilon} \vert >n \right\rbrace } c_1(x,t)dx dt +  \vert\Omega \vert^{\frac{1}{N^{\prime}}- \frac{1}{p^{\prime}} } \tilde{L}^{\frac{1}{p^{\prime}}} \right] \\
			& \quad \le 2 C^{\prime}(N,p) \left[ \frac{1}{C_1}  \int_0^{T} \int_{ \left\lbrace \vert u_{\epsilon} \vert >n \right\rbrace } c_1(x,t)dx dt +  \vert\Omega \vert^{\frac{1}{N^{\prime}}- \frac{1}{p^{\prime}} } \tilde{L}^{\frac{1}{p^{\prime}}} \right] 
		\end{aligned}
	\end{equation}
	For $\eta$ given in \eqref{3.23a}, we choose the integer $m$ to be large enough to guarantee that
	\begin{align}
		\frac{\eta C^{\prime}(N,p)}{C_1} \le \frac{1}{2}.
	\end{align}
	We denote $E_{n} = \left\lbrace  \vert u_{\epsilon} \vert >n\right\rbrace $. We observe that $\vert E_{n} \vert \rightarrow 0$ as $n \rightarrow + \infty$. Since $F_{\epsilon}$ strongly converges to $F$ in $L^{p^{\prime}}( 0,T; W^{-1,p^{\prime}}(\Omega))$, and $K \in L^1(Q)$, we obtain
	\begin{equation}\label{eq: 3.46}
		\begin{aligned}
				\lim_{n \rightarrow + \infty} \tilde{L} \rightarrow 0, \quad \text{as} \ \ n \to + \infty.
		\end{aligned}
	\end{equation}
	
	As a consequence, we can pass to the limit in \eqref{eq: 3.45} as $n$ tends to $+\infty$ to obtain
	\begin{equation}\label{eq: 3.47}
		\begin{aligned}
			\lim_{n \rightarrow + \infty} \| \vert \nabla R_{n}(u_{\epsilon})\vert^{p-1}\|_{L^1(0,T; L^{N^{\prime},\infty }(\Omega)) } =0. 
		\end{aligned}
	\end{equation}
	From \eqref{eq: 3.47}, it follows that
	\begin{equation}\label{eq: 3.48}
		\begin{aligned}
			\lim_{n \rightarrow + \infty} \tilde{M} \rightarrow 0, \quad \text{as} \ \ n \to + \infty.
		\end{aligned}
	\end{equation}
	
	In view of $(\ref{eq: 3.46})$ and $(\ref{eq: 3.48})$, the right hand side of $(\ref{eq: 3.43})$ tends to $0$ when $n \rightarrow \infty$. This concludes our proof.
	
\end{proof}
 
\begin{lemma}[Lemma~4.1, $\cite{Betta}$] \label{Le: 3.2} 
	Assume that $\Omega$ is an open subset of $\mathbb{R}^{N}$ with finite measure and that $1 <p<N$. Let $u$ be a measurable function satisfying $T_{k}(u) \in L^{p}(0,T,W^{1,p}(\Omega))$, for every positive $k$, and such that
	\begin{align*}
		\| \nabla T_{k}(u) \|_{L^{p}(Q) }^{p}  \le k M+L, \quad k >0,
	\end{align*}
	where $M$ and $L$ are given constants. Then $\vert \nabla u\vert^{p-1}\in L^{\infty}(0,T, L^{  N^{\prime}, \infty}(\Omega))$ and
	\begin{equation}\label{3.24}
		\begin{aligned}
			\int_0^{T}\| \vert \nabla u \vert^{p-1} \|_{L^{N^{\prime}, \infty}(\Omega) } dt  &\le C(N,p) \left[  M + \vert \Omega \vert^{\frac{1}{N^{\prime}}- \frac{1}{p^{\prime}} } L^{\frac{1}{p^{\prime}}}     \right].\\
		\end{aligned}
	\end{equation}
\end{lemma}

\begin{lemma}\label{Le: 3.1}
	For any fixed $k>0$,
	\begin{align}
		 T_{k}(u_{\epsilon}) \rightharpoonup T_{k}(u) \ \ \text{weakly in} \ \ L^{p}(0,T; W_0^{1,p}(\Omega)).
	\end{align}	
\end{lemma}
\begin{proof}

Inserting $T_{k}(u_{\epsilon})$ as a test function in \eqref{3.8}, and integrating in space and time on $Q$ we get
	\begin{equation}\label{3.18}
		\begin{aligned}
			&	\frac{1}{2} \| T_{k}(u_{\epsilon}) \|^2_{L^2(\Omega)} + \int_{Q} A_{\epsilon}(x,t,u_{\epsilon}, \nabla u_{\epsilon}) \cdot \nabla T_{k}(u_{\epsilon}) dx dt\\
			& \quad + \ \int_{Q} H_{\epsilon}(x,t,u_{\epsilon}, \nabla u_{\epsilon}) T_{k}(u_{\epsilon})dx dt     + \int_{Q} G_{\epsilon}(x,t,u_{\epsilon}) T_{k}(u_{\epsilon}) dxdt\\
			& = \frac{1}{2} \| T_{k}(u_{0\epsilon}) \|_{L^2(\Omega)}^2	 + \int_{Q} F_{\epsilon} \cdot \nabla T_{k}(u_{\epsilon}) dxdt
		\end{aligned}
	\end{equation}
	Assumption \eqref{3.10} on $A_{\epsilon}(x,t,u_{\epsilon},\nabla u_{\epsilon})$ leads to the estimate
	\begin{equation}\label{3.19}
		\begin{aligned}
			&	 \int_{Q} A_{\epsilon}(x,t,u_{\epsilon},\nabla u_{\epsilon}) \cdot \nabla T_{k}(u_{\epsilon}) dx dt\\
			&\quad = 	\int_0^{T}  \int_{\left\lbrace \vert u_{\epsilon} \vert \le k \right\rbrace } A_{\epsilon}(x,t,u_{\epsilon}, \nabla u_{\epsilon}) \cdot \nabla u_{\epsilon} dxdt\\
			& \quad \ge \int_0^{T} \int_{\left\lbrace \vert u_{\epsilon} \vert \le k \right\rbrace } \Big( \alpha\vert \nabla u_{\epsilon} \vert^{p} - \left( \vartheta_0(x,t) \vert u_{\epsilon} \vert\right)^{p} - K(x,t) \Big)  dxdt\\ 
			& \quad \ge \alpha   \int_{Q} \vert \nabla T_{k}(u_{\epsilon}) \vert^{p} dx dt- k^{p} \| \vartheta_0 \|^{p}_{L^{p}(Q)} - \|K\|_{L^1(Q)}.
		\end{aligned}
	\end{equation}	
		By \eqref{3.13}, assumption \eqref{3.7}, Young and H\"older inequalities we obtain the bound
	\begin{equation}\label{3.20}
		\begin{aligned}
			&	\left|  \int_{Q} H_{\epsilon}(x,t,u_{\epsilon}, \nabla u_{\epsilon}) T_{k}(u_{\epsilon}) dx dt\right| \\
			& \quad \le k  \int_{Q} \vert  H_{\epsilon}(x,t,u_{\epsilon}, \nabla u_{\epsilon}) \vert dx dt\\
			& \quad \le k  \int_0^{T} \int_{ \left\lbrace \vert u_{\epsilon} \vert \le k \right\rbrace } \vartheta_1 (x,t)\vert \nabla u_{\epsilon} \vert^{p-1} dxdt \\
			& \qquad+ k \int_0^{T}  \int_{ \left\lbrace \vert u_{\epsilon} \vert > k \right\rbrace } \vartheta_1(x,t) \vert \nabla u_{\epsilon} \vert^{p-1} dxdt + k\| c_1 \|_{L^1(Q)} \\
			& \quad \le \frac{ k^{p}2^{p/p^{\prime}}}{p \alpha^{p/p^{\prime}}} \|\vartheta_1\|^{p}_{L^{p}(Q)} + \frac{\alpha}{2p^{\prime}} \| \nabla T_{k}(u_{\epsilon}) \|^{p}_{L^{p}(Q)} \\
			& \qquad + k \| R_{m}\vartheta_1 \|_{ L^{\infty}(0,T; L^{N,1}(\Omega)) }  \int_0^{T}  \int_{ \left\lbrace \vert u_{\epsilon} \vert > k \right\rbrace } \vert \nabla u_{\epsilon} \vert^{p-1} dxdt\\
			& \qquad + k m  \int_0^{T}  \int_{ \left\lbrace \vert u_{\epsilon} \vert > k \right\rbrace } \vert \nabla u_{\epsilon} \vert^{p-1} dxdt + k \| c_1\|_{L^1(Q)}\\
			& \quad \le \frac{ k^{p}2^{p/p^{\prime}}}{p^{\prime} \alpha^{p/p^{\prime}}} \|\vartheta_1\|^{p}_{L^{p}(Q)} + \frac{\alpha}{2p^{\prime}} \| \nabla T_{k}(u_{\epsilon}) \|^{p}_{L^{p}(Q)} + k \| c_1\|_{L^1(Q)}\\
			& \qquad + k \eta \| \vert \nabla R_{k}(u_{\epsilon}) \vert^{p-1} \|_{ L^{1}(0,T;L^{N^{\prime}, \infty}(\Omega)) },
		\end{aligned}
	\end{equation}
	where $\eta$ is given by \eqref{3.23a}.

The sign condition \eqref{3.15} on $G_{\epsilon}$ yields
	\begin{equation}\label{3.21}
		\begin{aligned}
			\int_{Q} G_{\epsilon}(x,t,u_{\epsilon}, \nabla u_{\epsilon}) T_{k}(u_{\epsilon}) dxdt \ge 0.
		\end{aligned}
	\end{equation}
Finally, Young's inequality implies 
	\begin{equation}\label{3.22}
		\begin{aligned}
			& \int_{Q} F_{\epsilon} \cdot \nabla T_{k}(u_{\epsilon}) dxdt \\
			&\quad  \le \frac{\alpha}{2p} \| \nabla T_{k}(u_{\epsilon}) \|_{L^{p}(Q) }^{p}  + \frac{2^{p^{\prime}/p}}{p^{\prime} \alpha^{p^{\prime}/p}} \|F_{\epsilon}\|^{p^{\prime}}_{L^{p^{\prime}}(0,T; W^{-1, p^{\prime}}(\Omega))}.
		\end{aligned}
	\end{equation}
Combining \eqref{3.19}-\eqref{3.22}, we obtain
\begin{equation}\label{3.23}
	\begin{aligned}
	&	\frac{1}{2} \| T_{k}(u_{\epsilon}) \|^2_{L^2(\Omega)} +	\frac{\alpha}{2} \|  \nabla T_{k}(u_{\epsilon})\|_{L^{p}(Q)}^{p} \\
	& \quad \le \frac{1}{2} \| T_{k}(u_{0\epsilon}) \|_{L^2(\Omega)}^2 +k M_1 +L_1,
	\end{aligned}
\end{equation}
where $M_1$ and $L_1$ are defined by
\begin{equation*}
	\begin{aligned}
		&M_1 = \eta \| \vert \nabla R_{k}(u_{\epsilon}) \vert^{p-1} \|_{ L^1(0,T;L^{N^{\prime}, \infty}(\Omega)) } + \| c_1\|_{L^1(Q)},\\
		& L_1 = k^{p} \| \vartheta_0\|^{p}_{L^{p}(Q)} + \|K\|_{L^1(Q)} 
		 + \frac{ k^{p}2^{p/p^{\prime}}}{p \alpha^{p/p^{\prime}}} \|\vartheta_1\|^{p}_{L^{p}(Q)} + \frac{2^{p^{\prime}/p}}{p^{\prime} \alpha^{p^{\prime}/p}} \|F_{\epsilon}\|^{p^{\prime}}_{L^{p^{\prime}}(0,T; W^{-1, p^{\prime}}(\Omega))}.
	\end{aligned}
\end{equation*}
Using the estimate on $\| \vert \nabla R_{k}(u_{\epsilon}) \vert^{p-1} \|_{ L^1(0,T;L^{N^{\prime}, \infty}(\Omega)) } $ given in $(\ref{eq: 3.45})$, and by  \eqref{3.23}, it follows that $\nabla T_k(u_\epsilon)$ is bounded in $L^{p}(0,T; L^{p}(\Omega))$.  Following the argument in the proof of Theorem~2.1 in \cite{A. Porretta}, we can conclude the weak convergence of $T_k(u_\epsilon)$. 
	
\end{proof}


Previous a priori estimates allow us to extract a subsequence of $u_{\epsilon}$ (still indexed by $\epsilon$), as stated in the following Lemma.

\begin{lemma}\label{Le: 3.4}
There exists a subsequence (still indexed by $\epsilon$) of the sequence $u_{\epsilon}$ and an element $u$ belongs to $C^0([0,T], L^2(\Omega))$ such that 
\begin{align}
&	u_{\epsilon} \rightarrow u \quad \text{in} \quad C^0([0,T]; L^2(\Omega)); \label{3.42}\\ 
& T_{k}(u_{\epsilon}) \rightarrow T_{k}(u) \quad \text{in} \quad L^{p}(0,T; W_0^{1,p}(\Omega)); \label{3.43}\\
& H_{\epsilon} (x,t,u_{\epsilon}, \nabla u_{\epsilon}) \rightarrow H(x,t,u, \nabla u) \quad \text{strongly in} \quad L^1(Q); \label{3.44}\\
& G_{\epsilon}(x,t,u_{\epsilon}) \rightarrow G(x,t,u) \quad \text{strongly in} \quad L^1(Q); \label{3.45}
\end{align}
as $\epsilon$ tends to $0$, for any $k >0$.
\end{lemma}
We present here a brief proof of Lemma $\ref{Le: 3.4}$ for the sake of completeness.
\begin{proof}
 By considering the difference of equation \eqref{3.8} for $u_{\epsilon}$ with \eqref{3.8} for $u_{\nu}$  and using $\sigma^{-1} T_{\sigma}(u_{\epsilon}-u_{\nu})$ as a test function, as $\sigma\to 0$ we obtain the strong convergence of $u_{\epsilon}$ to $u$ in $C^0([0,T];L^2(\Omega))$, i.e. \eqref{3.42}.
 
By Lemma $\ref{Le: 3.1}$, $T_{k}(u_{\epsilon}) \rightarrow T_{k}(u)$ weakly in $L^{p}(0,T; W_0^{1,p}(\Omega))$ and by the asymptotic behavior of $T_{k}(u_{\epsilon})$ as $\epsilon$ tends to zero, we can deduce the strong convergence of $T_{k}(u_{\epsilon})$ to $T_{k}(u)$ in $L^{p}(0,T; W_0^{1,p}(\Omega))$, i.e. \eqref{3.43}.
 
By the definition of $H_{\epsilon}$ and \eqref{3.42}, we derive that
\begin{align}
	H_{\epsilon}(x,t,u_{\epsilon}, \nabla u_{\epsilon}) \rightarrow H(x,t,u,\nabla u) \quad \text{a.e. in} \quad Q.
\end{align}

Due to the growth condition \eqref{3.13} on $H_{\epsilon}$ and the generalized H\"older inequality, it follows that $H_{\epsilon}(x,t,u_{\epsilon}, \nabla u_{\epsilon})$ is equi-integrable. Thus, thanks to Vitali's lemma we deduce that
\begin{align}
	H_{\epsilon}(x,t,u_{\epsilon}, \nabla u_{\epsilon}) \rightarrow H(x,t,u,\nabla u) \quad \text{strongly} \quad L^1(Q).
\end{align}

Similarly, we can prove that 
\begin{align}
	G_{\epsilon}(x,t,u_{\epsilon}) \rightarrow G(x,t,u) \quad \text{strongly} \quad L^1(Q).
\end{align}
This concludes our proof.

\end{proof}

The next lemma establishes a monotonicity estimate that is crucial for the proof of the existence result.
\begin{lemma} \label{Le: 3.5}
	For any positive real number $k$, we have
	\begin{align}\label{3.49}
		 \lim_{\epsilon \rightarrow 0} \int_0^{T} \int_0^{t}\int_{\Omega} [  A(x,t, T_{k}(u_{\epsilon}), \nabla T_{k}(u_{\epsilon})) - A(x,t, T_{k}(u_{\epsilon}), \nabla T_{k}(u)    ]\cdot \nabla \left[ T_{k}(u_{\epsilon}) -T_{k}(u) \right] dxdsdt =0.
	\end{align}
\end{lemma}
\begin{proof}
For any $n \ge 1$, we consider an increasing sequence $S_{n}\in C^{\infty}(\mathbb{R})$ such that $S_{n}(r)=r$ for $r \le n$, $\mathrm{supp}S_{n}^{\prime} \subset [-(n+1), n+1]$ and $\| S_{n}^{\prime \prime}(r)\|_{L^{\infty}(\mathbb{R})} \le 1$. Pointwise multiplication of the approximating problem \eqref{3.8} by $S_{n}^{\prime}(u_{\epsilon})$ leads to 
	\begin{equation}\label{3.50}
		\begin{aligned}
			&(S_{n}(u_{\epsilon}))_{t}- \mathrm{div} \Big( A(x,t,u_{\epsilon}, \nabla u_{\epsilon}) S_{n}^{\prime}(u_{\epsilon}) \Big) + S_{n}^{\prime \prime}(u_{\epsilon}) A(x,t, u_{\epsilon}, \nabla u_{\epsilon})\cdot \nabla u_{\epsilon}\\
			& + H_{\epsilon}(x,t, u_{\epsilon}, \nabla u_{\epsilon}) S_{n}^{\prime}(u_{\epsilon}) + G_{\epsilon} (x,t,u_{\epsilon}) S_{n}^{\prime}(u_{\epsilon}) = - \mathrm{div} F_{\epsilon}(x,t) S_{n}^{\prime}(u_{\epsilon})
		\end{aligned}
	\end{equation}
	which is the renormalized formulation of \eqref{3.8} for $u_{\epsilon}$ with $S$ is replaced by $S_{n}$.
	
	Let $(T_{k}(u))_{\nu}$ be an approximation of $T_{k}(u)$. For any fixed $k>0$, we denote	
	\begin{align*}
		W_{\nu}^{\epsilon}= T_{k}(u_{\epsilon})- (T_{k}(u))_{\nu}.
	\end{align*}

We now make use of the test function $W_{\nu}^{\epsilon}$ in \eqref{3.50} and integrate over $(0,t)$ and then over $(0,T)$, we get
	\begin{equation}\label{3.51}
		\begin{aligned}
			&	\int_0^{T} \int_0^{t} \left\langle  \frac{\partial S_{n}(u_{\epsilon})}{\partial t}, W_{\nu}^{\epsilon}   \right\rangle dsdt - \int_0^{T} \int_0^{t} \int_{\Omega} S_{n}^{\prime}(u_{\epsilon}) A_{\epsilon}(x,t,u_{\epsilon}, \nabla u_{\epsilon}) \cdot \nabla W_{\nu}^{\epsilon} dxdsdt \\
			& \quad  +\int_0^{T} \int_0^{t} \int_{\Omega} S_{n}^{\prime \prime}(u_{\epsilon}) A_{\epsilon}(x,t,u_{\epsilon}, \nabla u_{\epsilon}) \cdot \nabla u_{\epsilon} W_{\nu}^{\epsilon} dxdsdt\\
			& \quad + \int_0^{T} \int_0^{t} \int_{\Omega} S_{n}^{\prime} (u_{\epsilon})H_{\epsilon} (x,t, u_{\epsilon}, \nabla u_{\epsilon}) W_{\nu}^{\epsilon} dxdsdt\\
			& \quad + \int_0^{T} \int_0^{t} \int_{\Omega} S_{n}^{\prime} (u_{\epsilon})G_{\epsilon} (x,t,u_{\epsilon}) W_{\nu}^{\epsilon} dxdsdt\\
			&= - \int_0^{T} \int_0^{t} \int_{\Omega} F_{\epsilon} (x,t)\cdot \nabla \Big( S_{n}^{\prime}(u_{\epsilon}) W_{\nu}^{\epsilon} \Big) dxdsdt.
		\end{aligned}
	\end{equation}
	
	To prove \eqref{3.49}, we need to prove 
	\begin{align}
		& \liminf_{\nu \rightarrow + \infty} \lim_{\epsilon \rightarrow 0} \int_0^{T} \int_0^{t} \left\langle  \frac{\partial S_{n}(u_{\epsilon})}{\partial t}, W_{\nu}^{\epsilon}   \right\rangle dsdt \ge 0 \quad \text{for any } n \ge k, \label{3.52}\\
		& \lim_{\nu \rightarrow + \infty} \lim_{\epsilon \rightarrow 0} \left| \int_0^{T} \int_0^{t} \int_{\Omega} S_{n}^{\prime \prime}(u_{\epsilon}) A_{\epsilon}(x,t, u_{\epsilon}, \nabla u_{\epsilon}) \cdot \nabla u_{\epsilon} W_{\nu}^{\epsilon} dxdsdt \right| =0, \label{3.53}\\
		& \lim_{\nu \rightarrow + \infty} \lim_{\epsilon \rightarrow 0} \int_0^{T} \int_0^{t} \int_{\Omega} S_{n}^{\prime} (u_{\epsilon})H_{\epsilon} (x,t,u_{\epsilon}, \nabla u_{\epsilon}) W_{\nu}^{\epsilon} dxdsdt =0, \label{3.54}\\
		& \lim_{ \nu \rightarrow + \infty} \lim_{\epsilon \rightarrow 0} \int_0^{T} \int_0^{t} \int_{\Omega} S_{n}^{\prime} (u_{\epsilon})G_{\epsilon} (x,t,u_{\epsilon}) W_{\nu}^{\epsilon} dxdsdt =0, \label{3.55}\\
		& \lim_{\nu \rightarrow + \infty} \lim_{\epsilon \rightarrow 0}   \int_0^{T} \int_0^{t} \int_{\Omega} F_{\epsilon} (x,t)\cdot \nabla \Big( S_{n}^{\prime}(u_{\epsilon}) W_{\nu}^{\epsilon} \Big) dxdsdt =0. \label{3.56}
	\end{align}
\textbf{Proof of \eqref{3.52}.} It follows immediately from Lemma~1 of $\cite{Blanchard}$.

\noindent\textbf{Proof of \eqref{3.53}.}
By the assumption \eqref{3.10} on $A_{\epsilon}$, $ \mathrm{supp}S_{n}^{\prime \prime} \subset [-(n+1), n+1] \cap [n, n+1]$, and  $W_{\nu}$ is bounded, we get
	\begin{equation*}
		\begin{aligned}
			&	\left|   \int_0^{T} \int_0^{t} \int_{\Omega} S_{n}^{\prime \prime}(u_{\epsilon}) A_{\epsilon}(x,t,u_{\epsilon}, \nabla u_{\epsilon}) \cdot \nabla u_{\epsilon} W_{\nu}^{\epsilon} dxdsdt \right| \\
			& \le T \| S_{n}^{\prime \prime} \|_{L^{\infty}(\mathbb{R})} \| W_{\nu}^{\epsilon} \|_{L^{\infty}(Q)} \int_{ \left\lbrace n \le \vert u_{\epsilon} \vert \le n+1 \right\rbrace } 	A_{\epsilon}(x,t,u_{\epsilon}, \nabla u_{\epsilon}) \cdot \nabla u_{\epsilon} dx dsdt\\
			& \le C \liminf_{\epsilon \rightarrow 0}	\int_{ \left\lbrace n \le \vert u_{\epsilon} \vert \le n+1 \right\rbrace } 	A_{\epsilon}(x,t,u_{\epsilon}, \nabla u_{\epsilon}) \cdot \nabla u_{\epsilon} dxdsdt\\
			& \le C^{\prime} \int_{ \left\lbrace n \le \vert u_{\epsilon} \vert \le n+1 \right\rbrace } \vert \nabla u_{\epsilon}\vert^{p} dxdsdt
		\end{aligned}
	\end{equation*}
 where $C^{\prime}$ is a constant independent of $n$. This can be achieved by using \eqref{3.27}-\eqref{3.29}. Lemma $\ref{Le: 3.4}$ allows us pass to the limit as $n$ tends to $+\infty$ to obtain \eqref{3.53}.

\noindent\textbf{Proof of \eqref{3.54}}. In view of \eqref{3.42}, \eqref{3.44} of Lemma \eqref{Le: 3.4}, we apply Lebesgue's convergence theorem to obtain
\begin{equation*}
	\begin{aligned}
	&	\lim_{\epsilon \rightarrow 0} \int_0^{T} \int_0^{t} \int_{\Omega} S_{n}^{\prime}(u_{\epsilon}) H_{\epsilon}(x,t,u_{\epsilon}, \nabla u_{\epsilon}) W_{\nu}^{\epsilon} dxdsdt \\
	& \quad = \int_0^{T} \int_0^{t} \int_{\Omega} S_{n}^{\prime}(u) H(x,t,u, \nabla u) (T_{k}(u) - (T_{k}(u))_{\nu}) dxdsdt.
	\end{aligned}
\end{equation*}
Since $(T_{k}(u))_{\nu} \rightarrow T_{k}(u)$ as $\nu\to +\infty$,  \eqref{3.54} holds true.


\noindent\textbf{Proof of \eqref{3.55}}. Using \eqref{3.42}, \eqref{3.45} of Lemma \ref{Le: 3.4} we get
\begin{equation*}
	\begin{aligned}
		&	\lim_{\epsilon \rightarrow 0} \int_0^{T} \int_0^{t} \int_{\Omega} S_{n}^{\prime}(u_{\epsilon}) G_{\epsilon}(x,t,u_{\epsilon}) W_{\nu}^{\epsilon} dxdsdt \\
		& \qquad = \int_0^{T} \int_0^{t} \int_{\Omega} S_{n}^{\prime}(u) G(x,t,u) (T_{k}(u) - (T_{k}(u))_{\nu}) dxdsdt
	\end{aligned}
\end{equation*}
Passing to the limit as $\nu\to+\infty$, we obtain \eqref{3.55}.

\noindent\textbf{Proof of \eqref{3.56}}. For any $n \ge 1$ and any $\nu>0$, we have
\begin{equation*}
	\begin{aligned}
		&	\lim_{\epsilon \rightarrow 0}  \int_0^{T} \int_0^{t} \int_{\Omega} F_{\epsilon}(x,t) \cdot \nabla \Big( S_{n}^{\prime}(u_{\epsilon}) W_{\nu}^{\epsilon} \Big) dxdsdt\\
		& \qquad =   \int_0^{T} \int_0^{t} \int_{\Omega} S_{n}^{\prime}(u) F(x,t) \Big( \nabla T_{k}(u) - \nabla (T_{k}(u))_{\nu} \Big) dxdsdt\\
		& \qquad \qquad +  \int_0^{T} \int_0^{t} \int_{\Omega} S_{n}^{\prime \prime}(u)F(x,t) \cdot \nabla T_{n+1}(u) \Big( T_{k}(u)- (T_{k}(u))_{\nu} \Big) dxdsdt,
	\end{aligned}
\end{equation*}
Since $ S_{n}^{\prime}(u) F \in L^{p^{\prime}}(Q)$, $ S_{n}^{\prime \prime}(u) F \cdot \nabla T_{n+1}(u) \in L^1(Q)$, and $(T_{k}(u))_{\nu}$ converges to $T_{k}(u) $ in $L^{p}(0,T, W_0^{1,p}(\Omega)) $ as $\nu \to + \infty$, we  obtain \eqref{3.56}.

\bigskip\noindent
Combining \eqref{3.52}-\eqref{3.56} we deduce that for any $k \ge 0$.
\begin{equation*}
	\begin{aligned}
		\lim_{n \rightarrow + \infty}  \limsup_{\nu \rightarrow  +\infty} \int_0^{T} \int_0^{t} \int_{\Omega} S_{n}^{\prime} (u_{\epsilon})A_{\epsilon}(x,t,u_{\epsilon}, \nabla u_{\epsilon}) \cdot \nabla \left[ T_{k}(u_{\epsilon})-T_{k}(u) \right]  dxdsdt \le 0.
	\end{aligned}
\end{equation*}

We observe that for any $k \le n$, 
\begin{equation}\label{3.57}
	\begin{aligned}
		&\int_0^{T} \int_0^{t} \int_{\Omega}  S_{n}^{\prime} (u_{\epsilon}) A_{\epsilon}(x,t,u_{\epsilon}, \nabla u_{\epsilon}) \cdot \nabla \left[  T_{k}(u_{\epsilon}) -  T_{k}(u) \right]  dxdsdt \\
		& \ge \int_0^{T} \int_0^{t}\int_{\Omega} A(x,t,T_{k}(u_{\epsilon}), \nabla T_{k}(u_{\epsilon})) \cdot \nabla \left[  T_{k}(u_{\epsilon}) -  T_{k}(u) \right]  dxdsdt \\
		& \quad -\int_0^{T}\int_{\Omega} S_{n}^{\prime}(u_{\epsilon}) A_{\epsilon}(x,t,u_{\epsilon}, \nabla u_{\epsilon})\cdot  \nabla T_{k}(u_{\epsilon})\chi_{\left\lbrace \vert u_{\epsilon}\vert >k \right\rbrace }dxdsdt
	\end{aligned}
\end{equation}
which yields
\begin{equation}\label{3.58}
	\begin{aligned}
		&\int_0^{T} \int_0^{t}\int_{\Omega}  S_{n}^{\prime}(u_{\epsilon}) A_{\epsilon}(x,t, u_{\epsilon}, \nabla u_{\epsilon}) \cdot \nabla \left[  T_{k}(u_{\epsilon})- T_{k}(u)  \right]  dxdsdt   \\
		& \ge \int_0^{T} \int_0^{t}\int_{\Omega} \left[  A(x,t,T_{k}(u_{\epsilon}), \nabla T_{k}(u_{\epsilon}))- A(x,t, T_{k}(u_{\epsilon}), \nabla T_{k}(u)) \right] \\
		& \qquad \qquad  \cdot  \nabla \left[  T_{k}(u_{\epsilon})-  T_{k}(u)  \right]  dxdsdt \\
		& \quad -\int_0^{T} \int_0^{t}\int_{\Omega} A(x,t,T_{k}(u_{\epsilon}), \nabla T_{k}(u_{\epsilon})) \cdot \nabla  \left[  T_{k}(u_{\epsilon})- T_{k}(u)  \right]  dxdsdt  \\
		& \quad+\int_0^{T} \int_0^{t}\int_{\Omega} A(x,t,T_{k}(u_{\epsilon}), \nabla T_{k}(u)) \cdot \nabla \left[   T_{k}(u_{\epsilon})-  T_{k}(u)  \right]    dxdsdt  \\
		& \quad-   \int_0^{T} \int_0^{t}\int_{\Omega}  S_{n}^{\prime}(u_{\epsilon}) A_{\epsilon}(x,t,u_{\epsilon}, \nabla u_{\epsilon}) \cdot \nabla T_{k}(u_{\epsilon})\chi_{\left\lbrace \vert u_{\epsilon}\vert >k \right\rbrace }dxdsdt.
	\end{aligned}
\end{equation}
Using the fact that $S_{n}^{\prime}$ has compact support and $\nabla T_{k}(u^{\epsilon})$ is bounded in $L^{p}(Q)$ for every $k>0$. Taking the limit in $\epsilon$ tends to zero, the last three terms of the right hand side of \eqref{3.58} goes to zero, thus \eqref{3.57} becomes
\begin{equation}\label{3.59}
	\begin{aligned}
		&\limsup_{\epsilon \rightarrow 0}  \int_0^{T} \int_0^{t}\int_{\Omega} \Big[ A(x,t,T_{k}(u_{\epsilon}),\nabla T_{k}(u_{\epsilon}))- A(x,t,T_{k}(u_{\epsilon}), \nabla T_{k}(u))\Big] \\
		& \qquad \qquad \cdot \nabla \left[ T_{k}(u_{\epsilon})-T_{k}(u) \right]  dxdsdt \\
		&\quad \le \limsup_{\epsilon \rightarrow  0}\int_0 ^{T} \int_0^{t}\int_{\Omega} S_{n}^{\prime}(u_{\epsilon})A_{\epsilon}(x,t,u_{\epsilon}, \nabla u_{\epsilon}) \nabla \left[ T_{k}(u_{\epsilon})-T_{k}(u)\right] dxdsdt.
	\end{aligned}
\end{equation}
Taking $\epsilon$ tends to $0$ in \eqref{3.59}, we obtain
\begin{equation}
	\begin{aligned}
			\limsup_{\epsilon \rightarrow 0} & \int_0^{T} \int_0^{t}\int_{\Omega} \Big[ A(x,t,T_{k}(u_{\epsilon}),\nabla T_{k}(u_{\epsilon}))- A(x,t,T_{k}(u_{\epsilon}),\nabla T_{k}(u)) \Big]\\
		& \qquad  \qquad \cdot  \nabla \left[ T_{k}(u_{\epsilon})-T_{k}(u)\right]  dxds dt=0
	\end{aligned}
\end{equation}
which concludes our proof.

\end{proof}

We will also need the following
\begin{lemma}\label{Le: 3.6}
	For fixed $k\ge 0$, as $\epsilon\to 0$ we have
\begin{align}\label{eq: 3.66}
	A(x,t,T_{k}(u_{\epsilon}), \nabla T_{k}(u_{\epsilon})) \cdot \nabla T_{k}(u_{\epsilon}) \rightharpoonup A(x,t,T_{k}(u), \nabla T_{k}(u)) \cdot \nabla T_{k}(u) \ \ \text{in } \quad L^1(Q).
\end{align}
\end{lemma}
\begin{proof}
	This result is shown in $\cite{Blanchard}$, so we omit the details of the proof.
\end{proof}

\subsection{Passing to the limit}

We are now in position to prove the existence result stated in Theorem \ref{Th: 3.1}. Let us show that the limiting solution $u$ satisfies \eqref{3.4}-\eqref{3.5}. 

Let $k$ be a positive real number and consider a function $S\in W^{2, \infty}(\mathbb{R})$ such that $S^{\prime}$ has compact support contained in $[-k,k]$. Pointwise multiplication of the approximating equation \eqref{3.8} with $S^{\prime}(u_{\epsilon})$ leads to
	\begin{equation}\label{3.63}
	\begin{aligned}
	&(S(u_{\epsilon}))_{t}- \mathrm{div} \Big( A_{\epsilon}(x,t,u_{\epsilon}, \nabla u_{\epsilon}) S^{\prime}(u_{\epsilon}) \Big) + S^{\prime \prime}(u_{\epsilon}) A_{\epsilon}(x,t,u_{\epsilon}, \nabla u_{\epsilon})\cdot \nabla u_{\epsilon}\\
		& +S^{\prime}(u_{\epsilon}) H_{\epsilon}(x,t,u_{\epsilon}, \nabla u_{\epsilon})  + S^{\prime}(u_{\epsilon}) G_{\epsilon} (x,t,u_{\epsilon})\\ 
	&	=  \mathrm{div} \Big( F_{\epsilon} (x,t)S_{n}^{\prime}(u_{\epsilon}) \Big) - S^{\prime \prime}(u_{\epsilon}) F_{\epsilon}(x,t) \cdot \nabla u_{\epsilon}.
	\end{aligned}
\end{equation}

We next pass to the limit as $\epsilon$ tends to $0$ in each term of  \eqref{3.63}.

By Lemma \ref{Le: 3.4}  $u_{\epsilon}$ converges strongly to $u$, moreover $S$ is bounded and continuous, so $S(u_{\epsilon})$ converges to $S(u)$ a.e in $Q$. Consequently, $(S_{n}(u_{\epsilon}))_{t}$ converges to $(S_{n}(u))_{t}$ as $\epsilon$ tends to $0$ in the distributional sense, specifically in the space ${\mathcal D}^{\prime}(Q)$. As supp$S^{\prime} \subset [-k,k]$, we have
\begin{equation*}
	\begin{aligned}
		&S^{\prime}(u_{\epsilon}) A_{\epsilon}(x,t,u_{\epsilon}, \nabla u_{\epsilon}) = 	S^{\prime}(u_{\epsilon}) A(x,t,T_{k}(u_{\epsilon}), \nabla T_{k}(u_{\epsilon})) \ \ \text{a.e. in } \ \ Q;\\
		& 	S^{\prime \prime}(u_{\epsilon}) A_{\epsilon}(x,t,u_{\epsilon}, \nabla u_{\epsilon})\cdot \nabla u_{\epsilon} = S^{\prime \prime}(u_{\epsilon}) A(x,t,T_{k}(u_{\epsilon}), \nabla T_{k}(u_{\epsilon})) \cdot \nabla T_{k}(u_{\epsilon}) \ \ \text{a.e. in} \ \ Q.
	\end{aligned}
\end{equation*}
Using again \eqref{3.42} and the fact that $S^{\prime}$ and $S^{\prime \prime}$ are both bounded, from Lemma \ref{Le: 3.6} we deduce that
\begin{equation}
	\begin{aligned}
			&S^{\prime}(u_{\epsilon}) A_{\epsilon}(x,t,u_{\epsilon}, \nabla u_{\epsilon}) \rightharpoonup  S^{\prime} (u)A(x,t,T_{k}(u), \nabla T_{k}(u)) \ \ \text{weakly in} \ \ L^{p^{\prime}}(Q);\\
			& S^{\prime \prime}(u_{\epsilon}) A_{\epsilon}(x,t,u_{\epsilon}, \nabla u_{\epsilon})\cdot \nabla u_{\epsilon}  \rightharpoonup S^{\prime \prime}(u) A(x,t,T_{k}(u), \nabla T_{k}(u)) \cdot \nabla T_{k}(u) \ \ \text{weakly in} \ \ L^1(Q).
	\end{aligned}
\end{equation} 
Note that the term $S^{\prime} (u) A(x,t,T_{k}(u), \nabla T_{k}(u))$ coincides here with $S^{\prime}(u) A(x,t,u, \nabla u)$, hence we identify $S^{\prime \prime}(u)A(x,t,T_{k}(u), \nabla T_{k}(u)) \cdot \nabla T_{k}(u) $ with the term $S^{\prime \prime}(u) A(x,t,u, \nabla u) \cdot \nabla u$ in  \eqref{3.4} of Definition \ref{Def: 3.1}.


By \eqref{3.42}, \eqref{3.44}, \eqref{3.45} and the bound on $S^{\prime}$ we have
\begin{equation}
	\begin{aligned}
	&	S^{\prime}(u_{\epsilon}) H_{\epsilon}(x,t,u_{\epsilon}, \nabla u_{\epsilon}) \rightarrow S^{\prime}(u) H(x,t,u, \nabla u) \ \ \text{ strongly in } \ \ L^1(Q);\\
	&S^{\prime}(u_{\epsilon}) G_{\epsilon}(x,t,u_{\epsilon}) \rightarrow S^{\prime}(u) G(x,t,u) \ \ \text{ strongly in } \ \ L^1(Q).
	\end{aligned}
\end{equation}


By \eqref{3.42} it follows that $F_{\epsilon}\to F$ and $S^{\prime}(u_{\epsilon})\to S^{\prime}(u)$ a.e. in $Q$ as $\epsilon\to 0$, thus
\begin{equation}
	\begin{aligned}
		&\mathrm{div}(F_{\epsilon}(x,t)S^{\prime}(u_{\epsilon})) \rightarrow \mathrm{div} ( F(x,t) S^{\prime }(u))  \ \ \text{strongly in}  \ \ L^{p^{\prime}}(0,T,W^{-1,p^{\prime}}(\Omega));\\
		&	S^{\prime \prime}(u_{\epsilon}) F_{\epsilon}(x,t)\cdot \nabla u_{\epsilon}= F_{\epsilon} \nabla S^{\prime}(u_{\epsilon}) \rightharpoonup F \cdot \nabla S^{\prime}(u)\ \ \text{weakly in} \ \ L^1(Q).
	\end{aligned}
\end{equation}

 Combining all of the above convergence results, we may pass to the limit as $\epsilon\to 0$ in \eqref{3.63} to conclude that the solution $u$ satisfies \eqref{3.4} in sense of Definition \ref{Def: 3.1}. 

The proof of \eqref{3.5} follows the same methodology as in $\cite{Murat-Blan})$, hence it is omitted. This completes the proof of Theorem \ref{Th: 3.1}.


\section{Uniqueness of renormalized solutions}
\begin{theorem}\label{Th: 4.1}
	Let us assume that conditions \eqref{1.2}-\eqref{1.7} are fulfilled, and that assumption \eqref{3.7} holds true. Then problem \eqref{1.1} has a unique renormalized solution.
\end{theorem}
The proof of Theorem \ref{Th: 4.1} will follow immediately from the following comparison lemma. 
\begin{lemma}\label{Le: 4.1}
	Assume that the initial data $u_{01},u_{02}\in L^2(\Omega)$ and the source terms  $F_1,F_2\in L^{p^{\prime}}(0,T; W^{-1, p^{\prime}}(\Omega))$. Assume furthermore that they satisfy (almost everywhere)
	\begin{align}
		& u_{01} \le u_{02},\\
		& F_1 \le F_2.
 	\end{align}
Then if $u_1$ and $u_2$ are two renormalized solution of problem \eqref{1.1} respectively for the data $(F_1,u_{01})$ and $(F_2, u_{02})$, we have
	\begin{align}\label{4.3}
		u_1 \le u_2  \quad \text{a.e. in } \  Q.
	\end{align}	
\end{lemma}

\begin{proof}
Our proof follows the ideas presented in the paper $\cite{Murat}$, providing a slight generalization to address the present context of noncoercive nonlinear parabolic problems with singularities and unbounded coefficients in the lower order terms.

For all $n \ge 2$, let $S_{n}$ be a smooth approximation of the truncation $T_{n}$ such that $\mathrm{supp} S_{n}^{\prime} \subset [-(n+1), n+1], \Vert S_{n}^{\prime} \Vert_{L^{\infty}(\mathbb{R})} \le \Vert h\Vert_{L^{\infty}(\mathbb{R})}$ and $S_{n}(r)=S_{n}(T_{n+1}(r))$. Here, for all $s$ in $\mathbb{R}$, the function $h$ is given by
\begin{align} h_n(s)=
	\begin{cases}
		1 \ \ \text{if} \ \ \vert s\vert \le n-1; \nonumber\\
		h(s-(n-1)\mathrm{sign}(s)), \ \ \text{if}\ \ \vert s \vert \ge n-1;
	\end{cases}
\end{align}	
where $sign(s)$ denotes the sign of $s$.

Let $C>0$ be fixed and consider $T_{C}^{+}(S_{n}(u_1)-S_{n}(u_2))=\max \{0, T_{C}(S_{n}(u_1)-S_{n}(u_2))\}$ as test function in the difference of equation \eqref{3.4} with respect to $u_1$ and $u_2$. Integrating over $(0,t)$ and then over $(0,T)$, we get		
\begin{equation}\label{4.4}
	\begin{aligned}
		 &\int_0^{T} \int_0^{t} \int_{\Omega}\left\langle (S_{n}(u_1))_{t} - (S_{n}(u_2))_{t}, T_{C}^{+}(S_{n}(u_1)-S_{n}(u_2))\right\rangle (s)dx ds dt \\
		 & \quad +I_{n}^{\left[ 1\right] }  + I_{n}^{\left[ 2\right] }  +I_{n}^{\left[ 3\right] }  +J_{n}^{\left[ 1\right] }  + J_{n}^{\left[ 2\right] }  \\
		&= \int_0^{T} \int_0^{t}\int_{\Omega} \left( F_1 S_{n}^{\prime}(u_1)- F_2 S_{n}^{\prime}(u_2)\right) \nabla T_{C}^{+}(S_{n}(u_1)-S_{n}(u_2)) dx ds dt 
	\end{aligned}
\end{equation}
here we denote $I_{n}^{ \left[ i\right] }, i =1,2,3$ and $J_{n}^{\left[ j\right] }, j =1, 2$ in the following
\begin{equation*}\label{eq: 4.5}
	\begin{aligned}
&	I_{n}^{\left[ 1\right] } =  \int_0^{T} \int_0^{t} \int_{\Omega} \left[ S_{n}^{\prime}(u_1) A(x,t,u_1, \nabla u_1) - S_{n}^{\prime}(u_2) A(x,t,u_2, \nabla u_2) \right] \cdot \nabla T_{C}^{+} (S_{n}(u_1) -S_{n}(u_2)) dxdsdt,\\	
& 	I_{n}^{\left[ 2\right] } =  \int_0^{T} \int_0^{t} \int_{\Omega} S_{n}^{\prime \prime}(u_1) A(x,t,u_1, \nabla u_1) \cdot \nabla u_1 T_{C}^{+} (S_{n}(u_1) -S_{n}(u_2)) dxdsdt,\\
& 	I_{n}^{\left[ 3\right] } =  \int_0^{T} \int_0^{t} \int_{\Omega} S_{n}^{\prime \prime}(u_2) A(x,t,u_2, \nabla u_2) \cdot \nabla u_2 T_{C}^{+} (S_{n}(u_1) -S_{n}(u_2))  dxdsdt,\\
& 	J_{n}^{\left[ 1\right] }= \int_0^{T} \int_0^{t} \int_{\Omega} \left[  S_{n}^{\prime}(u_1) H(x,t,u_1, \nabla u_1) - S_{n}^{\prime} (u_2)H(x,t,u_2, \nabla u_2) \right] T_{C}^{+} (S_{n}(u_1) -S_{n}(u_2)) dxdsdt,\\
&  	J_{n}^{\left[ 2\right] }=  \int_0^{T} \int_0^{t} \int_{\Omega} \left[  S_{n}^{\prime}(u_1) G(x,t,u_1)-  S_{n}^{\prime}(u_2) G(x,t,u_2)     \right] T_{C}^{+} (S_{n}(u_1) -S_{n}(u_2)) dxdsdt.	
	\end{aligned}
\end{equation*}
 For all $t\in [0,T]$, one has
\begin{align*}
	&\int_0^{t} \left\langle (S_{n}(u_1))_{t} - (S_{n}(u_2))_{t}, T_{C}^{+}(S_{n}(u_1)-S_{n}(u_2))\right\rangle (s) ds \nonumber \\
	&= \int_{\Omega}\tilde\phi_{C}([S_{n}(u_1)-S_{n}(u_2)]^{+}) (t)dx -\int_{\Omega} \tilde\phi_{C}([S_{n}(u_{01})-S_{n}(u_{02})]^{+})dx,
\end{align*}
where $\tilde\phi_C$ is defined in \eqref{tildephi}. Thus \eqref{4.4} becomes
\begin{equation} 
\begin{aligned}
	&\int_{Q}  \tilde\phi_{C}([S_{n}(u_1)-S_{n}(u_2)]^{+}) (t)dx dt\\ 
	& \quad +	I_{n}^{\left[ 1\right] } + 	I_{n}^{\left[ 2\right] } + 	I_{n}^{\left[ 3\right] }  + J_{n}^{\left[ 1\right] }+ 	J_{n}^{\left[ 2\right] } \\
	& \le  \int_0^{T} \int_0^{t}\int_{\Omega} \left( F_1 S_{n}^{\prime}(u_1)- F_2 S_{n}^{\prime}(u_2)\right)  \nabla T_{C}^{+}(S_{n}(u_1)-S_{n}(u_2)) dx ds dt \\
	& \quad +T\int_{\Omega}\tilde\phi_{C}([S_{n}(u_{01})-S_{n}(u_{02})]^{+})dx.
\end{aligned}		
\end{equation}
To prove \eqref{4.3}, we need to show
\begin{equation}\label{4.6}
	\begin{aligned}
	&	\lim_{n \rightarrow + \infty} I_{n}^{ \left[i \right] } =0, \quad i =1, 2,3 ;\\
	& 		\lim_{n \rightarrow + \infty} J_{n}^{ \left[j\right] } =0, \quad j =1, 2. 
	\end{aligned}
\end{equation}
The long technical details to prove \eqref{4.6} are postponed to the Appendix. From \eqref{4.6} we deduce that
\begin{equation}\label{4.7}
	\begin{aligned}
		&\int_{Q}  \tilde\phi_{C}([S_{n}(u_1)-S_{n}(u_2)]^{+}) (t)dx dt\\ 
		&\quad  \le  \int_0^{T} \int_0^{t}\int_{\Omega} \left( F_1 S_{n}^{\prime}(u_1)- F_2 S_{n}^{\prime}(u_2)\right)  \nabla T_{C}^{+}(S_{n}(u_1)-S_{n}(u_2)) dx ds dt \\
		& \qquad +T\int_{\Omega}\tilde\phi_{C}([S_{n}(u_{01})-S_{n}(u_{02})]^{+})dx.
	\end{aligned}		
\end{equation}
Since
\begin{equation}
	\begin{aligned}
		&S_{n}(u_1)-S_{n}(u_2) \rightarrow u_1 -u_2 \quad \text{ strongly in } \quad L^1(Q),\\
		& S_{n}(u_{01})-S_{n}(u_{02}) \rightarrow u_{01} -u_{02}  \quad \text{ strongly in } \quad L^2(\Omega),
	\end{aligned}
\end{equation}
it follows that
\begin{equation}
	\begin{aligned}
		\left( F_1 S_{n}^{\prime}(u_1)- F_2 S_{n}^{\prime}(u_2)\right) \cdot  \nabla T_{C}^{+}(S_{n}(u_1)-S_{n}(u_2)) \rightarrow (F_1 -F_2) \cdot  T_{C}^{+} (u_1 -u_2) 
	\end{aligned}
\end{equation}
strongly in $L^1(Q)$. Passing to the limit as $n\to+\infty$ in \eqref{4.7} we get
\begin{equation}\label{4.10}
	\begin{aligned}
		\int_{Q}  \tilde\phi_{C}([u_1- u_2]^{+}) (t)dx dt &\le   \int_0^{T} \int_0^{t}\int_{\Omega}(F_1 -F_2) \cdot  \nabla T_{C}^{+} (u_1 -u_2)  dxdsdt \\
		& \qquad + T\int_{\Omega}\tilde\phi_{C}([u_{01}- u_{02}]^{+})dx.
	\end{aligned}
\end{equation}
Since the right hand side of \eqref{4.10} is nonnegative, we have $\tilde\phi_{C}([u_1- u_2]^{+}) \le 0$. Since $C>0$ was arbitrary, this implies $u_1 \le u_2$ almost everywhere in $Q$. This concludes our proof.
\end{proof}

\appendix
\section{Proof of Lemma \ref{Le: 4.1} completed}

In this appendix we give a detailed proof of \eqref{4.6} to complete the proof of Lemma \ref{Le: 4.1}.

Let us consider first $I_{n}^{\left[1 \right] }$. We have
\begin{equation}
	\begin{aligned}
		I_{n}^{\left[ 1\right] }&= \int_{Q} (T-t)\chi_{\left\lbrace 0\le S_{n}(u_1)-S_{n}(u_2) \le C\right\rbrace }\left[S_{n}^{\prime}(u_1) A(x,t,u_1, \nabla u_1) - S_{n}^{\prime}(u_2) A(x,t,u_2, \nabla u_2)\right] \\
		& \qquad	\times \nabla \left[ S_{n}(u_1)-S_{n}(u_2)\right] dxdt,\\
		& = I_{n}^{\left[1,1 \right] } + I_{n}^{\left[ 1,2\right] } + I_{n}^{\left[ 1,3\right] },
	\end{aligned}	
\end{equation}
where 
\begin{equation}
	\begin{aligned}
			I_{n}^{\left[ 1,1\right] }&= \int_{Q} (T-t)\chi_{\left\lbrace 0\le S_{n}(u_1)-S_{n}(u_2) \le C\right\rbrace }[A(x,t,u_1, \nabla S_{n}(u_1))- A(x,t,u_2, \nabla S_{n}(u_2))] \\
			& \qquad \times \nabla \left[ S_{n}(u_1)-S_{n}(u_2) \right] dxdt,	
	\end{aligned}
\end{equation}

\begin{equation}
	\begin{aligned}
			I_{n}^{\left[ 1,2\right] } &=  \int_{Q} (T-t)\chi_{\left\lbrace 0\le S_{n}(u_1)-S_{n}(u_2) \le C\right\rbrace }[S_{n}^{\prime}(u_1)A(x,t,u_1, \nabla u_1)- A(x,t,u_1, \nabla S_{n}(u_1))] \\
			&  \qquad \times  \nabla \left[  S_{n}(u_1)- S_{n}(u_2) \right]  dxdt,
	\end{aligned}
\end{equation}

\begin{equation}
	\begin{aligned}
		I_{n}^{\left[ 1,3\right] } &= - \int_{Q} (T-t)\chi_{\left\lbrace 0\le S_{n}(u_1)-S_{n}(u_2) \le C\right\rbrace }[S_{n}^{\prime}(u_1)A(x,t,u_2, \nabla u_2)- A(x,t,u_2, \nabla S_{n}(u_2))] \\
		&  \qquad \times  \nabla \left[  S_{n}(u_1)- S_{n}(u_2) \right]  dxdt.
	\end{aligned}
\end{equation}
By the assumption \eqref{1.4} on $A$, we obtain
\begin{align}
		I_{n}^{\left[ 1,1\right] } \geq 0.
\end{align}
Since $S$ is smooth, then
\begin{align*}
	S_{n}^{\prime}(u_1) A(x,t,u_1, \nabla u_1) = A(x,t,u_1, \nabla S_{n}(u_1))
\end{align*}
almost everywhere except on the subset of $\left\lbrace (x,t): n \le \vert u_1 \vert \le n+1 \right\rbrace $. 

As $\vert S_n \vert _{L^\infty(\mathbb{R})} \leq \vert h \vert_{L^\infty(\mathbb{R})}$ and  $\mathrm{supp} S_n \subset [-(n+1), n+1]$ we can deduce that for any $n > C$
\begin{align*}
	&\chi_{\left\lbrace 0\le S_{n}(u_1)-S_{n}(u_2) \le C\right\rbrace }\chi_{\left\lbrace n\le \vert u_1 \vert \le n+1\right\rbrace }S_{n}^{\prime}(u_2)\nonumber\\
	& =\chi_{\left\lbrace 0\le S_{n}(u_1)-S_{n}(u_2) \le C\right\rbrace } \chi_{\left\lbrace n\le \vert u_1 \vert \le n+1\right\rbrace }\chi_{\left\lbrace n-C\le \vert u_2\vert \le n+1  \right\rbrace } S_{n}^{\prime}(u_2).
\end{align*}	


  The boundedness assumption \eqref{1.3} on $A$ leads to
\begin{equation}\label{eq: 4.15}
	\begin{aligned}
			I_{n}^{\left[ 1,2\right] } &= \int_0^{T} \int_{ \left\lbrace n \le \vert u_1 \vert \le n+1 \right\rbrace }(T-t)\chi_{\left\lbrace 0\le S_{n}(u_1)-S_{n}(u_2) \le C\right\rbrace } \\
			& \quad \times [S_{n}^{\prime}(u_1)A(x,t,u_1, \nabla u_1)- A(x,t,u_1, \nabla S_{n}(u_1))] \times  \nabla \left[  S_{n}(u_1)- S_{n}(u_2) \right]  dxdt  \\
			& \le T  \int_0^{T} \int_{ \left\lbrace n \le \vert u_1 \vert \le n+1 \right\rbrace }  \left[  S_{n}^{\prime}(u_1)A(x,t,u_1, \nabla u_1)- A(x,t,u_1, \nabla S_{n}(u_1)) \right]   \\
			& \qquad \times \left[    \| h \|_{L^{\infty}(\mathbb{R})} \vert \nabla u_1 \vert -  \chi_{\left\lbrace n-C\le \vert u_2\vert \le n+1  \right\rbrace }     \| h \|_{L^{\infty}(\mathbb{R})} \vert \nabla u_2 \vert   \right]  dxdt \\
			& \le \hat{I}_{n}^{\left[1,2 \right] } + \tilde{I}_{n}^{\left[1,2\right] }
	\end{aligned}
\end{equation}
where
\begin{equation*}
	\begin{aligned}
		 \hat{I}_{n}^{\left[1,2 \right] } & = T  \| h \|_{L^{\infty}(\mathbb{R})}  \int_0^{T} \int_{ \left\lbrace n \le \vert u_1 \vert \le n+1 \right\rbrace } \left[  \beta \vert \nabla u_1 \vert^{p-1} + (b \vert u_1 \vert)^{p-1} + a_0 \right] \\
		 & \qquad \times \left[   \| h \|_{L^{\infty}(\mathbb{R})} \vert \nabla u_1 \vert -  \chi_{\left\lbrace n-C\le \vert u_2\vert \le n+1  \right\rbrace } \| h \|_{L^{\infty}(\mathbb{R})} \vert \nabla u_2 \vert  \right] dxdt ,
	\end{aligned}
\end{equation*}
\begin{equation*}
	\begin{aligned}
		\tilde{I}_{n}^{\left[1,2\right] } & = T \int_0^{T} \int_{ \left\lbrace n \le \vert u_1 \vert \le n+1 \right\rbrace }  \left[  \beta    \| h \|_{L^{\infty}(\mathbb{R})}^{p-1} \vert \nabla u_1 \vert^{p-1} + (b \vert u_1 \vert)^{p-1} + a_0  \right] \\
		& \qquad \times \left[   \| h \|_{L^{\infty}(\mathbb{R})} \vert \nabla u_1 \vert -\chi_{\left\lbrace n-C\le \vert u_2\vert \le n+1  \right\rbrace }  \| h \|_{L^{\infty}(\mathbb{R})} \vert \nabla u_2 \vert  \right] dxdt .
	\end{aligned}
\end{equation*}
To simplify notation, let us set
\begin{align*}
	\chi_{n,C} = \chi_{\left\lbrace n-C \le \vert u_2 \vert \le n+1\right\rbrace }.
\end{align*}

We first consider $ \hat{I}_{n}^{\left[1,2\right] }$. We have
\begin{equation}\label{eq: 4.18}
	\begin{aligned}
		\hat{I}_{n}^{\left[1,2\right] } & =  T  \| h \|_{L^{\infty}(\mathbb{R})} \int_0^{T} \int_{ \left\lbrace n \le \vert u_1 \vert \le n+1 \right\rbrace }  \left[  \beta \vert \nabla u_1 \vert^{p-1} + (b \vert u_1 \vert)^{p-1} + a_0 \right] \\
		& \qquad \times \left[  \| h \|_{L^{\infty}(\mathbb{R})} \vert \nabla u_1 \vert -\chi_{n,C} \| h \|_{L^{\infty}(\mathbb{R})} \vert \nabla u_2 \vert   \right]  dxdt\\
		 & = II_1 + II_2 +II_3,
	\end{aligned}
\end{equation}
where
\begin{equation*}
	\begin{aligned}
		& II_1 = T\beta \| h \|_{L^{\infty}(\mathbb{R})} \int_0^{T} \int_{ \left\lbrace n \le \vert u_1 \vert \le n+1 \right\rbrace }  \vert \nabla u_1 \vert^{p-1}  \times \left[  \| h \|_{L^{\infty}(\mathbb{R})} \vert \nabla u_1 \vert - \chi_{n,C}  \| h \|_{L^{\infty}(\mathbb{R})} \vert \nabla u_2 \vert   \right]  dxdt ,\\ 
		& II_2 = T \| h \|_{L^{\infty}(\mathbb{R})}  \int_0^{T} \int_{ \left\lbrace n \le \vert u_1 \vert \le n+1 \right\rbrace }    (b \vert u_1 \vert)^{p-1} \times \left[   \| h \|_{L^{\infty}(\mathbb{R})} \vert \nabla u_1 \vert - \chi_{n,C} \| h \|_{L^{\infty}(\mathbb{R})} \vert \nabla u_2 \vert  \right]  dxdt ,\\ 
		& II_3 = T  \| h \|_{L^{\infty}(\mathbb{R})}  \int_0^{T} \int_{ \left\lbrace n \le \vert u_1 \vert \le n+1 \right\rbrace } a_0 
		\times \left[   \| h \|_{L^{\infty}(\mathbb{R})} \vert \nabla u_1 \vert - \chi_{n,C} \| h \|_{L^{\infty}(\mathbb{R})} \vert \nabla u_2 \vert  \right] dxdt.
	\end{aligned}
\end{equation*}
 By Young inequality we get
\begin{equation}
	\begin{aligned}
		II_1 &=  T \beta \| h \|_{L^{\infty}(\mathbb{R})} \int_0^{T} \int_{ \left\lbrace n \le \vert u_1 \vert \le n+1 \right\rbrace }    \vert \nabla u_1 \vert^{p-1}  \\
		& \quad \times \left[  \| h \|_{L^{\infty}(\mathbb{R})} \vert \nabla u_1 \vert - \chi_{n,C}  \| h \|_{L^{\infty}(\mathbb{R})} \vert \nabla u_2 \vert  \right] dxdt\\ 
		& \le  T \| h \|_{L^{\infty}(\mathbb{R})}^2  \beta \int_0^{T} \int_{ \left\lbrace n \le \vert u_1 \vert \le n+1 \right\rbrace } \vert \nabla u_1 \vert^{p} dxdt \\
		& \quad +  T \| h \|_{L^{\infty}(\mathbb{R})}^2  \beta \frac{1}{p^{\prime}} \int_0^{T} \int_{ \left\lbrace n \le \vert u_1 \vert \le n+1 \right\rbrace } \vert \nabla u_1 \vert^{p} dxdt\\
		 & \quad +  T \| h \|_{L^{\infty}(\mathbb{R})}^2  \beta   \frac{1}{p} \int_0^{T} \int_{ \left\lbrace n \le \vert u_2\vert \le n+1 \right\rbrace } \vert \nabla u_2 \vert^{p} dxdt  \\
		 & \le  T \| h \|_{L^{\infty}(\mathbb{R})}^2  \beta \left( 1 + \frac{1}{p^{\prime}} \right)   \int_0^{T} \int_{ \left\lbrace n \le \vert u_1 \vert \le n+1 \right\rbrace } \vert \nabla u_1 \vert^{p} dxdt\\
		 & \quad + T \| h \|_{L^{\infty}(\mathbb{R})}^2  \beta   \frac{1}{p} \int_0^{T} \int_{ \left\lbrace n-C \le \vert u_2 \vert \le n+1 \right\rbrace } \vert \nabla u_2 \vert^{p} dxdt.
	\end{aligned}
\end{equation}

For the term $II_2$, we have
\begin{equation}\label{eq: 4.22}
	\begin{aligned}
		II_2 & = T \| h \|_{L^{\infty}(\mathbb{R})}  \int_0^{T} \int_{ \left\lbrace n \le \vert u_1 \vert \le n+1 \right\rbrace }     \Big(  b \vert u_1 \Big)^{p-1} \times \left[  \| h \|_{L^{\infty}(\mathbb{R})} \vert \nabla u_1 \vert - \chi_{n,C} \| h \|_{L^{\infty}(\mathbb{R})} \vert \nabla u_2 \vert  \right]  dxdt \\ 
		&\le  T\| h\|_{L^{\infty}(\mathbb{R})}  \int_0^{T} \int_{ \left\lbrace n \le \vert u_1 \vert \le n+1 \right\rbrace } \left[ \Big( R_{m}b \vert u_1 \vert \Big)^{p-1} + (T_{m} b \vert u_1 \vert)^{p-1}\right] \\
		& \qquad \times \left[   \| h \|_{L^{\infty}(\mathbb{R})} \vert \nabla u_1 \vert - \chi_{n,C}  \| h \|_{L^{\infty}(\mathbb{R})} \vert \nabla u_2 \vert  \right]  dxdt\\ 
		&  \le II_2^{\left[ 1\right] } +  II_2^{\left[ 2\right] },
	\end{aligned}
\end{equation}
where we denote
\begin{equation*}
	\begin{aligned}
		II_2^{\left[ 1\right] } &= T\| h\|_{L^{\infty}(\mathbb{R})}  \int_0^{T} \int_{ \left\lbrace n \le \vert u_1 \vert \le n+1 \right\rbrace } \Big( R_{m}b \vert u_1 \vert \Big)^{p-1} \\
		& \qquad \times \left[   \| h \|_{L^{\infty}(\mathbb{R})} \vert \nabla u_1 \vert - \chi_{n,C}  \| h \|_{L^{\infty}(\mathbb{R})} \vert \nabla u_2 \vert   \right]  dxdt,\\ 
	\end{aligned}
\end{equation*}
\begin{equation*}
	\begin{aligned}
			II_2^{\left[ 2\right] } & = T\| h\|_{L^{\infty}(\mathbb{R})}  \int_0^{T} \int_{ \left\lbrace n \le \vert u_1 \vert \le n+1 \right\rbrace }  (T_{m} b \vert u_1 \vert)^{p-1} \\
			& \qquad \times \left[   \| h \|_{L^{\infty}(\mathbb{R})} \vert \nabla u_1 \vert - \chi_{n,C}  \| h \|_{L^{\infty}(\mathbb{R})} \vert \nabla u_2 \vert   \right] dxdt.
	\end{aligned}
\end{equation*}
By the generalized H\"older inequality, Young inequality and the Sobolev embedding theorem, the first term of the right hand side of $(\ref{eq: 4.22})$ can be estimated as follows:
 \begin{equation}\label{eq: 4.25}
 	\begin{aligned}
 		&	II_2^{\left[ 1\right] } =T\| h\|_{L^{\infty}(\mathbb{R})}  \int_0^{T} \int_{ \left\lbrace n \le \vert u_1 \vert \le n+1 \right\rbrace } \left[ \Big( R_{m}b \vert u_1 \vert \Big)^{p-1} \right]\\
 		&  \qquad \times \left[   \| h \|_{L^{\infty}(\mathbb{R})} \vert \nabla u_1 \vert - \chi_{n,C}  \| h \|_{L^{\infty}(\mathbb{R})} \vert \nabla u_2 \vert   \right]  dxdt\\
 		& \quad   \le T\| h\|_{L^{\infty}(\mathbb{R})}^2  \int_0^{T} \int_{ \left\lbrace n \le \vert u_1 \vert \le n+1 \right\rbrace } (R_{m}b)^{p-1} \vert u_1 \vert^{p-1} \vert \nabla u_1 \vert dxdt\\
 		& \qquad + T\| h\|_{L^{\infty}(\mathbb{R})}^2  \int_0^{T} \int_{ \left\lbrace n \le \vert u_1 \vert \le n+1 \right\rbrace } \chi_{n,C}\Big( R_{m}b \vert u_1 \vert \Big)^{p-1} \vert \nabla u_2 \vert dxdt\\
 		& \quad \le T \| h\|_{L^{\infty}(\mathbb{R})}^2  \int_0^{T} \chi_{\left\lbrace n \le \vert u_1 \vert \le n+1 \right\rbrace } \| R_{m}b \|_{L^{N, \infty}(\Omega) }^{p-1} \| u_1 \|_{ p^{*},p}^{p-1} \| \nabla u_1 \|_{p} dt \\
 		& \qquad + T\| h\|_{L^{\infty}(\mathbb{R})}^2  \frac{1}{p^{\prime}} \int_0^{T} \int_{ \left\lbrace n \le \vert u_1 \vert \le n+1 \right\rbrace } \Big( R_{m}b \vert u_1 \vert \Big)^{p} dxdt \\
 	  &\qquad 	+ T\| h\|_{L^{\infty}(\mathbb{R})}^2 \frac{1}{p} \int_0^{T} \int_{ \left\lbrace n-C \le \vert u_2\vert \le n+1 \right\rbrace } \vert \nabla u_2 \vert^{p} dxdt\\
 	  & \quad \le T \| h\|_{L^{\infty}(\mathbb{R})}^2 S_{N,p}^{p-1}\| R_{m}b \|_{L^{\infty}(0,T; L^{ N, \infty}(\Omega))} ^{p-1} \int_0^{T} \int_{ \left\lbrace n \le \vert u_1 \vert \le n+1 \right\rbrace } \vert \nabla u_1 \vert^{p} dxdt \\
 	  & \qquad + T \| h\|_{L^{\infty}(\mathbb{R})}^2 \frac{1}{ p^{\prime} } S_{N,p}^{p} \| R_{m}b \|_{L^{\infty}(0,T; L^{N, \infty}(\Omega))}^{p} \int_0^{T} \int_{ \left\lbrace n \le \vert u_1 \vert \le n+1 \right\rbrace } \vert \nabla u_1 \vert^{p} dxdt\\
 	  & \qquad + T\| h\|_{L^{\infty}(\mathbb{R})}^2 \frac{1}{p} \int_0^{T} \int_{ \left\lbrace n-C \le \vert u_2\vert \le n+1 \right\rbrace } \vert \nabla u_2 \vert^{p} dxdt\\
 	  & \quad \le 2T \| h\|_{L^{\infty}(\mathbb{R})}^2 \frac{\alpha}{2p} \left( 1+ \frac{1}{p^{\prime}} \right)   \int_0^{T} \int_{ \left\lbrace n \le \vert u_1 \vert \le n+1 \right\rbrace } \vert \nabla u_1 \vert^{p} dxdt \\
 	 & \qquad + T\| h\|_{L^{\infty}(\mathbb{R})}^2 \frac{1}{p} \int_0^{T} \int_{ \left\lbrace n-C \le \vert u_2\vert \le n+1 \right\rbrace } \vert \nabla u_2 \vert^{p} dxdt.
 	\end{aligned}
 \end{equation}

Similarly, for the second term $II_2^{\left[ 2\right] } $, we have
\begin{equation}\label{eq: 4.26}
	\begin{aligned}
	&	II_2^{\left[ 2\right] } = T \| h\|_{L^{\infty}(\mathbb{R})}  \int_0^{T} \int_{ \left\lbrace n \le \vert u_1 \vert \le n+1 \right\rbrace } (T_{m} b \vert u_1 \vert)^{p-1} \\
		& \qquad \times \Big(  \| h \|_{L^{\infty}(\mathbb{R})} \vert \nabla u_1 \vert - \chi_{n,C}  \| h \|_{L^{\infty}(\mathbb{R})} \vert \nabla u_2 \vert     \Big) dxdt\\
		& \quad \le T  \| h\|_{L^{\infty}(\mathbb{R})}^2 m^{p-1}  \int_0^{T} \int_{ \left\lbrace n \le \vert u_1 \vert \le n+1 \right\rbrace } \vert u_1 \vert^{p-1} \vert \nabla u_1 \vert dxdt \\
		& \qquad + T  \| h\|_{L^{\infty}(\mathbb{R})}^2 m^{p-1} \int_0^{T} \int_{ \left\lbrace n \le \vert u_1 \vert \le n+1 \right\rbrace }\chi_{n, C} \vert u_1 \vert^{p-1} \vert \nabla u_2 \vert dxdt\\
		& \quad \le T   \| h\|_{L^{\infty}(\mathbb{R})}^2 m^{p-1} \| 1 \|_{L^{N, \infty}(\Omega)  }^{p-1} S_{N,p}^{p-1} \int_0^{T} \int_{ \left\lbrace n \le \vert u_1 \vert \le n+1 \right\rbrace } \vert \nabla u_1 \vert^{p} dxdt\\
		& \qquad + T  \| h\|_{L^{\infty}(\mathbb{R})}^2 m^{p-1} \frac{1}{p^{\prime}} \int_0^{T} \int_{ \left\lbrace n \le \vert u_1 \vert \le n+1 \right\rbrace } \vert u_1 \vert^{p} dxdt \\
		& \qquad +  T  \| h\|_{L^{\infty}(\mathbb{R})}^2 m^{p-1} \frac{1}{p}  \int_0^{T} \int_{ \left\lbrace n-C \le \vert u_2 \vert \le n+1 \right\rbrace } \vert \nabla u_2 \vert^{p} dxdt\\
		& \quad \le T   \| h\|_{L^{\infty}(\mathbb{R})}^2 m^{p-1} \vert \Omega \vert^{\frac{p-1}{N}}  \left(1+ \frac{1}{p^{\prime}} \right)  S_{N,p}^{p-1} \int_0^{T} \int_{ \left\lbrace n \le \vert u_1 \vert \le n+1 \right\rbrace } \vert \nabla u_1 \vert^{p} dxdt\\
		&  \qquad +  T  \| h\|_{L^{\infty}(\mathbb{R})}^2 m^{p-1} \frac{1}{p}  \int_0^{T} \int_{ \left\lbrace n-C \le \vert u_2 \vert \le n+1 \right\rbrace } \vert \nabla u_2 \vert^{p} dxdt.
	\end{aligned}
\end{equation}
Taking to the limit as $n$ tends to $+\infty$, and using Lemma $\ref{Le: 3.4}$, we get
\begin{align}
	& \lim_{n \rightarrow + \infty} II_1 =0,\\
	&\lim_{n \rightarrow + \infty} II_2 =  \lim_{n \rightarrow + \infty} 	II_2^{\left[ 1\right] } + \lim_{n \rightarrow + \infty} 	II_2^{\left[ 2\right] } =   0,\\
	& \lim_{n \rightarrow + \infty} II_3 = 0.
\end{align}
As a consequence, we obtain
\begin{align}\label{eq: 4.29}
	\lim_{n \rightarrow + \infty} \hat{I}_{n}^{\left[ 1,2\right] } = 0.
\end{align}
Proceeding in the same way for $ \hat{I}_{n}^{\left[1,2 \right] }$, and pass to the limit as $n \rightarrow + \infty$, we obtain
\begin{align}\label{eq: 4.30}
	\lim_{n \rightarrow + \infty} 	\tilde{I}_{n}^{\left[1,2 \right] }  =0.
\end{align}
From $(\ref{eq: 4.15})$, using $(\ref{eq: 4.29})-(\ref{eq: 4.30})$, we conclude that
\begin{align}
	\lim_{n \rightarrow + \infty} I_{n}^{ \left[ 1,2\right] } =0.
\end{align}
By interchanging the roles of $u_1$ and $u_2$ and following a similar procedure as for $I_{n}^{[1,2]}$, we can deduce that $\lim_{n\rightarrow\infty} I_{n}^{[1,3]} = 0$. Together with the limits of $I_{n}^{[1,1]}$ and $I_{n}^{[1,2]}$, we can establish that 
\begin{equation}\label{eq: 4.32}
	\begin{aligned}
		\lim_{n\rightarrow\infty} I_{n}^{[1]} = 0.
	\end{aligned}
\end{equation}


Let us now provide an estimate of $I_{n}^{ \left[2\right] }$. As before, due to the ellipticity condition $(\ref{1.3})$ on $A$, by the Sobolev embedding theorem, we get
\begin{equation}\label{eq: 4.33}
	\begin{aligned}
	&	I_{n}^{\left[2 \right] } = \int_{Q} (T-t) S_{n}^{\prime \prime} A(x,t,u_1, \nabla u_1)\cdot \nabla u_1 T_{C}^{+}(S_{n}(u_1)-S_{n}(u_2)) dxdt\\
		& \quad \le T C  \| h^{\prime}\|_{L^{\infty}(\mathbb{R})} \int_0^{T} \int_{ \left\lbrace n \le \vert u_1 \vert \le n+1 \right\rbrace }  \left[ \alpha \vert \nabla u_1 \vert^{p} - \Big( b \vert u_1 \vert \Big)^{p} - K(x,t) \right] dxdt\\
		& \quad \le T C M \| h^{\prime}\|_{L^{\infty}(\mathbb{R})} \int_0^{T} \int_{ \left\lbrace n \le \vert u_1 \vert \le n+1 \right\rbrace }   \vert \nabla u_1 \vert^{p} dxdt \\
		& \qquad + T C  \| h^{\prime}\|_{L^{\infty}(\mathbb{R})}  \int_0^{T} \int_{ \left\lbrace \vert u_1 \vert >n\right\rbrace  } K(x,t)dxdt
	\end{aligned}
\end{equation}
with $M$ is defined by
\begin{align*}
	M = \frac{\alpha}{2p^{\prime}}- \vert \Omega \vert^{\frac{p}{N}} S_{N,p}^{p}.
\end{align*}
By Lemma $\ref{Le: 3.4}$, it is possible to pass to the limit in $(\ref{eq: 4.33})$ to obtain
\begin{align}\label{eq: 4.34}
	\lim_{n \rightarrow \infty} I_{n}^{\left[2 \right] }  =0.
\end{align}

Let us now consider the term $J_{n}^{\left[ 1\right] }$. By the growth assumption \eqref{1.5} on $H$, and by the generalized H\"older inequality, we deduce that
\begin{equation}\label{eq: 4.35}
	\begin{aligned}
  J_{n}^{\left[ 1\right] }&= \int_0^{T} \int_0^{t} \int_{\Omega} \left[  S_{n}^{\prime}(u_1) H(x,t,u_1, \nabla u_1) - S_{n}^{\prime} (u_2)H(x,t,u_2, \nabla u_2) \right] \\
  & \qquad \times T_{C}^{+} (S_{n}(u_1) -S_{n}(u_2)) dxdsdt\\
  & \quad \le TC   \| h\|_{L^{\infty}(\mathbb{R})}   \int_0^{T}\int_{ \left\lbrace \vert u_1 \vert > n \right\rbrace }\Big( c_0 \cdot \vert \nabla u_1 \vert^{p-1} + c_1\Big) dxdt\\
  & \qquad + TC   \| h\|_{L^{\infty}(\mathbb{R})}  \int_0^{T}\int_{ \left\lbrace \vert u_2 \vert > n \right\rbrace }  \Big( c_0  \cdot \vert \nabla u_2 \vert^{p-1} + c_1\Big) dxdt\\
  & \quad \le 2TC \| h\|_{L^{\infty}(\mathbb{R})} \int_{ \left\lbrace \vert u_1 \vert > n \right\rbrace }  c_1 dxdt+  J_{n}^{\left[ 1,1\right] } +  J_{n}^{\left[ 1,2\right] }, \\
	\end{aligned}
\end{equation}
where
\begin{equation}
	\begin{aligned}
		&J_{n}^{\left[ 1,1\right] } = TC \| h\|_{L^{\infty}(\mathbb{R})} \int_0^{T} \int_{ \left\lbrace \vert u_1 \vert > n \right\rbrace }   c_0 \cdot \vert \nabla u_1 \vert^{p-1}  dxdt,\\
		& J_{n}^{\left[ 1,2\right] } = T C \| h\|_{L^{\infty}(\mathbb{R})} \int_0^{T} \int_{ \left\lbrace \vert u_2 \vert > n \right\rbrace }   c_0 \cdot \vert \nabla u_2 \vert^{p-1}  dxdt.
	\end{aligned}
\end{equation}

 By the generalized H\"older inequality, we get
\begin{equation}
	\begin{aligned}
		J_{n}^{\left[1,1 \right] } &= T C \| h\|_{L^{\infty}(\mathbb{R})} \int_0^{T}\int_{ \left\lbrace \vert u_1 \vert > n \right\rbrace }
		c_0 \cdot \vert \nabla u_1 \vert^{p-1}  dxdt\\
		& \le TC \| h\|_{L^{\infty}(\mathbb{R})} \int_0^{T}\int_{ \left\lbrace \vert u_1 \vert > n \right\rbrace }    R_{m}c_0 \cdot \vert \nabla u_1 \vert^{p-1} dxdt\\
		& \quad+ T C \| h\|_{L^{\infty}(\mathbb{R})} \int_0^{T}\int_{ \left\lbrace \vert u_1 \vert > n \right\rbrace }   T_{m}c_0 \cdot \vert \nabla u_1 \vert^{p-1} dxdt\\
		& \le TC\| h\|_{L^{\infty}(\mathbb{R})} \eta \int_0^{T}  \| \vert \nabla R_{n}u_1 \vert^{p-1} \|_{L^{N^{\prime}, \infty}(\Omega)}dt\\
		& \le C_{\star} \int_0^{T} \| \vert \nabla R_{n}u_1 \vert^{p-1} \|_{L^{N^{\prime}, \infty}(\Omega)}dt,
	\end{aligned}
\end{equation}
where $C_{\star}= TC\| h\|_{L^{\infty}(\mathbb{R})} \eta$ with $\eta$ is given by \eqref{3.23a}. \\
Thanks to $(\ref{eq: 3.47})$, and passing to the limit as $n \rightarrow + \infty$, we can derive that
\begin{equation}\label{eq: 4.38}
	\begin{aligned}
	&	\lim_{n \rightarrow \infty} J_{n}^{\left[1,1 \right] } =0 
	\end{aligned}
\end{equation}
By employing a similar approach for $J_{n}^{\left[1,2 \right] } $, while interchanging the roles of $u_1$ and $u_2$, we get
\begin{equation}\label{eq: 4.39}
	\begin{aligned}
		&	\lim_{n \rightarrow \infty} J_{n}^{\left[1,2 \right] } =0 
	\end{aligned}
\end{equation}
From $(\ref{eq: 4.35})$, using $(\ref{eq: 4.38})$ and $(\ref{eq: 4.39})$, we deduce that
\begin{equation}\label{eq: 4.40}
	\begin{aligned}
		&	\lim_{n \rightarrow \infty} J_{n}^{\left[1 \right] } =0.
	\end{aligned}
\end{equation}	
Finally, the sign condition \eqref{1.7} on $G$ implies that
\begin{equation}\label{eq: 4.41}
	\begin{aligned}
		 	\lim_{n \rightarrow \infty} J_{n}^{\left[2 \right] } =0. 	
	\end{aligned}
\end{equation}	
This concludes the proof. \qed

\end{document}